\documentclass[nonatbib]{ipart}

\makeatletter
\let\c@author\relax
\makeatother

\RequirePackage{hyperref}
\usepackage{graphicx}
\usepackage{tikz}
\usepackage{nicefrac}
\usetikzlibrary{arrows}
\usetikzlibrary{decorations.markings}
\usepackage{subfig}

\newcommand{\R}{{I\!\!R}}
\newcommand{\X}{{\bf x}}
\newcommand{\Y}{{\bf y}}
\newcommand{\V}{{\bf v}}
\newcommand{\set}[1]{\left\{ #1 \right\}}
\newcommand{\size}[1]{\left| #1 \right|}
\newcommand{\norm}[1]{\left\| #1 \right\|}
\newcommand{\N}{\mathbf{N}}
\newcommand{\real}{\mathbf{R}}

\newcommand{\mySub}[1]{\paragraph*{\it{\textbf{#1.}}}}

\usepackage{biblatex}
\theoremstyle{plain}
\newtheorem{theorem}{Theorem}[section]

\newtheorem{lemma}[theorem]{Lemma}

\newtheorem*{comp_problem}{Computational Problem}

\theoremstyle{definition}
\newtheorem{definition}[theorem]{Definition}
\newtheorem{problem}{Open Problem}[section]

\theoremstyle{remark}
\newtheorem{example}[theorem]{Example}

\usepackage{suffix}

\usepackage{graphicx} 
\usepackage{xcolor}

\pubyear{2023}
\volume{0}
\issue{0}
\firstpage{1}
\lastpage{43}

\begin{document}

\begin{frontmatter}
\title[Seven open problems in applied combinatorics]{Seven open problems in applied combinatorics}
\thankstext{T1}{Information Release PNNL-SA-182814.
}

\begin{aug}
    \author{\fnms{Sinan G.} \snm{Aksoy}\ead[label=e1]{sinan.aksoy@pnnl.gov}},
    \address{Pacific Northwest National Laboratory\\
             \printead{e1}}
    \author{\fnms{Ryan} \snm{Bennink}\ead[label=e2]{benninkrs@ornl.gov}},
    \address{Oak Ridge National Laboratory\\
             \printead{e2}}
    \author{\fnms{Yuzhou} \snm{Chen}\ead[label=e13]{yuzhou.chen@temple.edu}},
    \address{Temple University\\
             \printead{e13}}
    \author{\fnms{Jos\'{e}} \snm{Fr\'{i}as}\ead[label=e11]{frias4@cimat.mx}},
    \address{University of Texas at Dallas\\
             \printead{e11}}
    \author{\fnms{Yulia R.} \snm{Gel}\ead[label=e8]{ygl@utdallas.edu}},
    \address{University of Texas at Dallas\\
             \printead{e8}}
    \author{\fnms{Bill} \snm{Kay}\ead[label=e3]{william.kay@pnnl.gov}},
    \address{Pacific Northwest National Laboratory\\
             \printead{e3}}
    \author{\fnms{Uwe} \snm{Naumann}\ead[label=e4]{naumann@stce.rwth-aachen.de}},
    \address{RWTH Aachen University\\
             \printead{e4}\\
             }
    \author{\fnms{Carlos} \snm{Ortiz Marrero}\ead[label=e10]{carlos.ortizmarrero@pnnl.gov}},
    \address{Pacific Northwest National Laboratory\\
             \printead{e10}}
    \author{\fnms{Anthony V.} \snm{Petyuk}\ead[label=e9]{anthony.petyuk@gmail.com}},
    \address{Hanford High School\\
             \printead{e9}}
    \author{\fnms{Sandip} \snm{Roy}\ead[label=e7]{sandip@wsu.edu}},
    \address{Washington State University\\
             \printead{e7}\\
             }
    \author{\fnms{Ignacio} \snm{Segovia-Dominguez}\ead[label=e12]{Ignacio.SegoviaDominguez@UTDallas.edu}},
    \address{University of Texas at Dallas\\
             Jet Propulsion Laboratory, Caltech\\
             \printead{e12}}
    \author{\fnms{Nate} \snm{Veldt}\ead[label=e5]{nveldt@tamu.edu}},
    \address{Texas A\&M University\\
             \printead{e5}\\
             }
    \and
    \author{\fnms{Stephen J.} \snm{Young}
            \ead[label=e6]{stephen.young@pnnl.gov}%
            }
    \address{Pacific Northwest National Laboratory\\
             \printead{e6}\\
             }
\end{aug}
\received{\sday{3} \smonth{1} \syear{2023}}

\begin{abstract}
We present and discuss seven different open problems in applied combinatorics. The application areas relevant to this compilation include quantum computing, algorithmic differentiation, topological data analysis, iterative methods, hypergraph cut algorithms, and power systems. 
\end{abstract}




\end{frontmatter}

\begin{refsection}
While the development of combinatorial mathematics has always been spurred by applied problems, in recent years both the scope and specificity of applied combinatorics has increased dramatically. In concrete and varied ways, combinatorial approaches are tackling problems across the sciences: from using graph clustering to design reduced order models and controls for complex fluid flows~\cite{meena2018network,meena2021identifying}, employing Ramanujan graphs as interconnection topologies of linkages between elements of supercomputers \cite{aksoy2021ramanujan, young2022spectralfly}, predicting drug interactions using graph neural networks (GNNs)~\cite{knutson2022decoding}, to detecting anomalies in cybersecurity data using graph and hypergraph centrality measures \cite{aksoy2021directional, joslyn2020hypergraph}. 

In celebration and promotion of this eclecticism, here we compile seven open problems in applied combinatorics posed by researchers in academia and government. The problems are organized in self-contained sections, authored by the attributed submitters. Each consists of a problem statement, discussion of relevant application areas, and any partial progress or prior work. The problems are as follows:

\begin{itemize}
    \item[\ref{sec:frias}.] {\bf The Dowker Complex in Metric Graphs}: Fr\'{i}as, Segovia-Dominguez, Chen, and Gel propose problems in topological data analysis (TDA) related to Dowker complexes. 
    \item[\ref{sec:kay}.]{\bf An Application of Probabilistic Combinatorics to Quantum Circuit Expressiveness}: Kay and Bennink introduce a problem in probabilistic combinatorics relevant to quantum circuits. 
    \item[\ref{sec:veldt}.] {\bf The Computational Complexity of the 4-uniform Hypergraph Minimum $s$-$t$ Cut Problem}: Veldt outlines a problem on the computational complexity of hypergraph minimum $s$-$t$ cuts.
    \item[\ref{sec:naumann1}-\ref{sec:naumann2}.] {\bf The Edge and Vertex Elimination Problems in Directed Acyclic Graphs} and {\bf Data Flow Reversal Problems}: Naumann poses several problems involving directed acyclic graphs, motivated by algorithmic differentiation of numerical programs.
    \item[\ref{sec:ratfish}.] {\bf Price of Asynchrony}: Ortiz Marrero and Young pose problems to study how asynchronous updates affect the convergence rate of iterative methods for large-scale systems in scientific computing. 
\item [\ref{sec:aksoy}.] {\bf $(n-k)$-contingent Zero-Forcing for Power Grids}: Aksoy, Petyuk, Roy, and Young propose problems on a variant of graph zero-forcing relevant to structural stability in power grids.
\end{itemize}

We have curated this sample of open problems with several criteria in mind. First, we have chosen problems that are timely and new, in the sense that they are motivated by recent research and not likely to be well-known by the larger community. Second, while far from comprehensive, our varied selection aims to illustrate the diversity of areas engaged by applied combinatorics. Third, while much of combinatorics plausibly ``has applications”, we have searched for problems having clear, specific impact to fields outside pure mathematics, including computing, industry, data science, and scientific software. We thank the reader for their interest.

\printbibliography[heading=subbibliography]
\end{refsection}

\begin{refsection}

\section{The Dowker Complex in Metric Graphs}\label{sec:frias}
\begin{flushright}
{\it Jos\'{e} Fr\'{i}as, Ignacio Segovia-Dominguez, Yuzhou Chen, Yulia R.\ Gel }
\end{flushright}

\mySub{Background}

Topological data analysis (TDA) is a set of methods in computational topology that continues to receive increasing attention in data sciences. One of the most relevant tools in TDA is persistent homology (PH)~\cite{CEH2006,ELZ2000, ZC2005}. The aim of PH is to obtain topological information of a finite metric space using a \emph{filtered simplicial complex} $\mathcal{K}$ with vertices in the data set, namely, there exists a sequence of simplicial complexes $\mathcal{K}_{0}\subset\mathcal{K}_{1}\subset \cdots\subset \mathcal{K}_{n-1}\subset\mathcal{K}_{n}=\mathcal{K}$. A non-trivial element $\gamma\in H_{p}(\mathcal{K}_{i})$, in the $p$-dimensional homology group of $\mathcal{K}_{i}$, is commonly called a \emph{$p$-dimensional topological feature} and is associated with a point $(i_b,j_d) \in \mathbb{R}^2$,  where $i_b$ and $j_d$ are the \emph{birth} and \emph{death} of $\gamma$, respectively (i.e., the indices of the simplicial complex at which $\gamma$ is first and last observed, respectively, that is, when it becomes trivial in homology or is merged to another topological feature that was born before). The lifespan of the topological feature is $j_d-i_b$.
The \emph{persistence diagram} $PD$ of the filtered simplicial complex $\mathcal{K}$ is the multiset $PD(\mathcal{K})=\{(i_b,j_d) \in \mathbb{R}^2\mid i_b<j_d\}$, where the multiplicity of $(i_b,j_d)$ is the number of topological features in the filtered simplicial complex $\mathcal{K}$ that are born and die at $i_b$ and $j_d$, respectively. A persistence diagram summarizes the evolution of the homology groups of $\mathcal{K}$ along the filtration. Two persistence diagrams can be compared using metrics such as the bottleneck or the Wasserstein distances (please see~\cite{CDO2014,CEH2006, ELZ2000, ZC2005} for a detailed exposition on persistence diagrams and distances in PH). 

Metric spaces induced by weighted graphs are particularly interesting due to their combinatorial foundations and  their applications in data analysis~\cite{ABB2020,LFWX2022}. A weighted graph $\mathcal{G}=(\mathcal{V},\mathcal{E},\omega)$ contains a set of vertices $\mathcal{V}$, an edge set $\mathcal{E}\subset \mathcal{V}\times \mathcal{V}$, and a weight function $\omega:\mathcal{E}\rightarrow \mathbb{R}_{+}$. We say that a graph $\mathcal{G}$ is \emph{simple} if it does not contain  self-edges nor multiple edges. Given a path $\gamma$ in a weighted graph $\mathcal{G}$, the \emph{length} of $\gamma$ is the sum of the weights of the edges in $\gamma$.  A connected weighted graph $\mathcal{G}$ is then endowed with the \emph{geodesic distance} $d_{\mathcal{G}}:\mathcal{ V}\times \mathcal{V} \rightarrow \mathbb{R}_{\geq 0}$, defined on a pair of vertices $u,v\in \mathcal{V}$ as the minimum length among all the paths connecting $u$ and $v$. The set $(\mathcal{V}, d_{\mathcal{G}})$ is a finite metric space. \par 

 One of the most widely used simplicial complexes in TDA is the Vietoris-Rips simplicial complex~\cite{CDO2014, ZC2005}. One important property of this complex is that it is completely determined by its $1$-skeleton, which makes it suitable for computations. 
 
\begin{definition}[{Vietoris-Rips Complex}]\label{df:vrc}
Let $\mathcal{G}=(\mathcal{V},\mathcal{E}, \omega)$ be a weighted simple graph with induced geodesic distance $d_{\mathcal{G}}$. For $\alpha \in \mathbb{R}_{\geq 0}$, we define the \textbf{Vietoris-Rips complex} $VR_{\alpha}(\mathcal{G})$ as the abstract simplicial complex with vertices in $\mathcal{V}$ and, for $k\geq 2$, a $k$-simplex  $\sigma=[x_{0},x_{1},\ldots,x_{k}]\in VR_{\alpha}(\mathcal{G})$ if and only if $d_{\mathcal{G}}(x_{i},x_{j})\leq \alpha$ for $0\leq i\leq j\leq k$.
\end{definition}

However, if the cardinality of $\mathcal{V}$ is large, the number of simplices in $VR_{\alpha}(\mathcal{G})$,  for large values $\alpha$, could be excessively large to be analyzed using computational tools. An alternative to deal with the scalability problem is to construct another simplicial complex with set of vertices $L\subset \mathcal{V}$, where $|L|<|\mathcal{V}|$. One of such simplicial complexes is the witness complex~\cite{ABB2020,BGO2009, DC2004, DFW2015}. 

\begin{definition}[{Witness Complex}]\label{df:wc}
Let $\mathcal{G}=(\mathcal{V},\mathcal{E}, \omega)$ be a weighted simple graph with induced geodesic distance $d_{\mathcal{G}}$, and take subsets $L,W\subset \mathcal{V}$. For $\alpha \in \mathbb{R}_{\geq 0}$, let $Wit_{\alpha}(W,L)$ be the abstract simplicial complex with set of vertices $L$, and a simplex $\sigma\in Wit_{\alpha}(W,L)$ if and only if for every $\tau\subseteq \sigma$ there exists $w\in W$ such that $d_{\mathcal{G}}(w,l)\leq d_{\mathcal{G}}(w,l') +\alpha$ for all $l\in\tau$ and $l'\in L\setminus\tau$.  The complex $Wit_{\alpha}(W,L)$ is called the \textbf{witness complex} of $\mathcal{G}$ with set of \textbf{witnesses} $W$ and \textbf{landmarks} $L$.
\end{definition}

Given a sequence of non-decreasing values $0\leq\alpha_0\leq\alpha_1\leq\cdots \leq\alpha_{n-1}\leq\alpha_n$, the associated Vietoris-Rips complexes constructed on $\mathcal{G}$ satisfy $VR_{\alpha_{0}}(\mathcal{G})\subseteq VR_{\alpha_{1}}(\mathcal{G}),\subseteq \cdots \subseteq VR_{\alpha_{n}}(\mathcal{G})$, namely, the sequence of scale values define a filtered simplicial complex. If we select sets of landmarks and witnesses $L,W\subset\mathcal{V}$, there exists also a filtration  $Wit_{\alpha_{0}}(L,W)\subseteq Wit_{\alpha_{1}}(L,W)\subseteq \cdots \subseteq Wit_{\alpha_{n}}(L,W)$. Then persistent homology can be computed in both filtered simplicial complexes. However, computation of witness complexes on graphs is challenging due to the difficulty to determine a witness of some simplexes, as the next illustrative example shows.

\begin{figure}
  \centering
    \includegraphics[width=10cm]{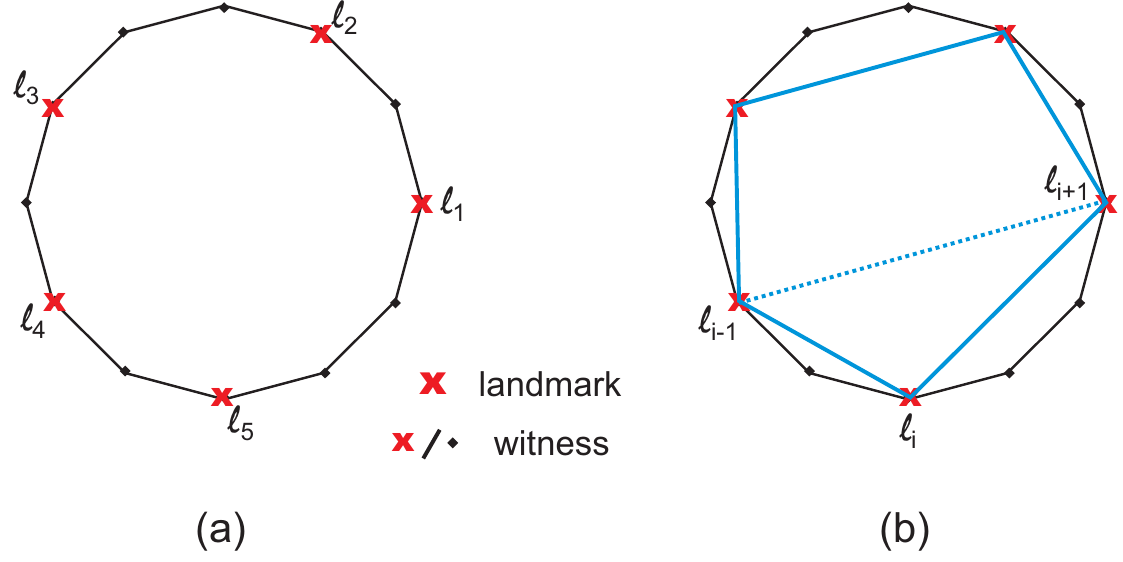}
  \caption{Landmarks in a cycle graph.}
  \label{fig:f1}
\end{figure}

\begin{example}
Suppose we are interested in the $1$-dimensional persistent homology of the $1$-weighted cycle graph $\mathcal{C}_{m}$, $m\geq 4$, with vertices set $\mathcal{V}$. Let $L\subset\mathcal{V}$ be a set of landmarks and $W=\mathcal{V}$ be a set of witnesses (see Figure~\ref{fig:f1}a). It is not difficult to determine the smallest scale value $\alpha> 0$ such that $Wit_{\alpha}(L,W)$ contains the $1$-simplex  $\{l_{i},l_{i+1}\}$, for any  two consecutive landmarks $l_{i},l_{i+1}\in L$ (in Figure~\ref{fig:f1}b the $1$-simplexes are represented by edges connecting consecutive landmarks). These $1$-simplexes form a $1$-dimensional topological feature $\gamma$ that is born at $\alpha$. However, given a sequence of three consecutive landmarks $l_{i-1},l_{i},l_{i+1}\in L$, it is combinatorially more difficult to determine a witness and a scale value at which the $1$-simplex $\{l_{i-1},l_{i+1}\}$ (dotted edge in Figure~\ref{fig:f1}b), appears in the filtered witness complex, or, even more difficult, when $\gamma$ vanishes.  
\end{example}

The Dowker complex~\cite{CDO2014, CM2018}  is a construction obtained after relaxing the definition of the witness complex, Definition~\ref{df:wc}, to make it easier to compute:   

\begin{definition}[{Dowker Complex}]\label{df:dc}
Let $\mathcal{G}=(\mathcal{V},\mathcal{E}, \omega)$ be a weighted simple graph with induced geodesic distance $d_{\mathcal{G}}$, and let $L,W\subset \mathcal{V}$.
 For $\alpha \in \mathbb{R}_{\geq 0}$, let $Dow_{\alpha}(W,L)$ be the abstract simplicial complex with set of vertices $L$ and a simplex $\sigma\in Dow_{\alpha}(W,L)$ if and only if there exists $w\in W$ such that $d_{\mathcal{G}}(w,l)\leq \alpha$ for all $l\in\sigma$. The simplicial complex $Dow_{\alpha}(W,L)$ is called the \textbf{Dowker complex} of $\mathcal{G}$ with sets of witnesses and landmarks $W$ and $L$, respectively.
\end{definition}

Note that witness and Dowker complexes constructions strongly depend on the selection of the sets $L$ and $W$. In~\cite{DC2004}, two algorithms to define the set of landmarks $L$ are developed: \emph{random} and \emph{maxmin} algorithms. More recently,
\cite{ABB2020} proposes the $\epsilon$-nets algorithm  to obtain a set of landmarks $L\subset \mathcal{V}$ satisfying both properties of sparsity and proximity to the set of points $\mathcal{V}\setminus L$, thereby allowing for a trade-off between computational complexity and potential information loss. 

\begin{definition}[{$\epsilon$-net}]\label{df:en}
 Let $(\mathcal{V},d_{\mathcal{G}})$ be the finite metric space obtained from the weighted graph $\mathcal{G}=(\mathcal{V},\mathcal{E},\omega)$. Given $L=\{u_{1},u_{2},\ldots,u_{l}\}\subset \mathcal{V}$ and $\epsilon \geq 0$ then:
\begin{itemize}
    \item [(i)] The set $L$ is an \textbf{$\epsilon$-sample} of $\mathcal{G}$ if the collection $\{\mathcal{N}(u_{i})\}_{i=1}^{l}$ of closed $\epsilon$-neighborhoods of points in $L$ covers $\mathcal{V}$, i.e. for any $v\in \mathcal{V}$ there exists $u_{j}\in L$ such that $d_{\mathcal{G}}(v,u_{j})\leq \epsilon$.
    \item[(ii)] $L$ is \textbf{$\epsilon$-sparse} if for any two distinct points $u_{i},u_{j}\in L$, their distance $d_{\mathcal{G}}(u_{i},u_{j})>\epsilon$.
    \item[(iii)] The set $L$ is an \textbf{$\epsilon$-net} of $\mathcal{G}$ if it is an $\epsilon$-sample of $\mathcal{G}$ and is $\epsilon$-sparse.
\end{itemize}
\end{definition}

\mySub{Open Problems}

Along the present section, $\mathcal{G}=(\mathcal{V},\mathcal{E},\omega)$ is a weighted simple graph and $d_{\mathcal{G}}$ is the geodesic distance induce by $\mathcal{G}$ on $\mathcal{V}$. The construction of Dowker simplicial complexes on the metric space $(\mathcal{V},d_{\mathcal{G}})$ relies on the selection of the sets of landmarks and witnesses. Usually the set of witnesses is taken to be the set of vertices, $W=\mathcal{V}$, or the complement of landmarks, $W=\mathcal{V}\setminus L$. 
Furthermore, there are numerous possibilities of selecting a set of landmarks.
The existence of an $\epsilon$-net defining the set of landmarks guarantees some combinatorial and geometric properties, as shown in~\cite{ABB2020}. Furthermore, \cite{ABB2020} proved that for a weighted graph $\mathcal{G}=(\mathcal{V},\mathcal{E},\omega)$ and a particular $\epsilon$, there exists an $\epsilon$-net for $(\mathcal{V},d_{\mathcal{G}})$ whose cardinality admits a bound depending on $\epsilon$, the diameter of the graph and the number of vertices (Theorem~3). However, this bound does not pretend to be any close to a sharp bound, while its further analysis has a number of important implication for graph learning applications in machine learning and statistics. As a result, we formulate the first open problem on the Dowker complexes on graphs as follows:

\begin{problem}\label{pr:1}
Given a weighted simple graph  $\mathcal{G}=(\mathcal{V},\mathcal{E},\omega)$ and $\epsilon>0$, determine  optimal upper and lower bounds for the number of elements of an $\epsilon$-net for $\mathcal{G}$. 
\end{problem} 

Three different methods to obtain $\epsilon$-nets in metric graphs, namely, the greedy, iterative and \textit{SPTpruning} algorithms, are proposed in~\cite{ABB2020}. These algorithms are experimentally compared with respect to the number of landmarks they produce for the same graph, including their execution time. 
As such, the next natural question arises:

\begin{problem}\label{pr:2}
For a weighted simple graph  $\mathcal{G}=(\mathcal{V},\mathcal{E},\omega)$, propose an algorithm to construct an $\epsilon$-net for $\mathcal{G}$ such that:
\begin{itemize}
\item[(i)]  It is based on a lower number of landmarks than that of the existing algorithms;
\item[(ii)] The time of execution of the algorithm is lower than the existing algorithms.
\end{itemize}
 \end{problem} 

Since the real-world datasets are usually noisy and incomplete, a suitable choice of landmarks in the Dowker complexes has multi-fold potential benefits in applications. 
First, it tends to improve performance of graph learning algorithms by focusing on the most essential
topological characteristics (i.e. the ``data skeleton"). Second, it helps enhance robustness to noise/data perturbations.
Third, by reducing computational costs, it opens a perspective of utilizing TDA for downstream tasks in graph learning. Indeed, as computing persistence homology on large graphs is often infeasible due to prohibitive computational costs,
the Dowker complex construction offers a promising alternative by the means of the suitable landmark selection. Additionally, the data separation in landmarks and witnesses trigger an intrinsic spatial decomposition that match with the natural local divisions needed when dealing with modern massive parallel computing~\cite{LEWIS2015}.

Given the current proliferation of the Vietoris-Rips complex in TDA, the natural question to address is
the comparison of persistent homology obtained from the Vietoris-Rips  and  Dowker complexes constructions on the same weighted graph $\mathcal{G}$. We denote with $PD^{i}(\mathcal{K})$ the persistence diagram of dimension $i$ for the filtered simplicial complex $\mathcal{K}$. In the case of metric graphs,  $1$-dimensional homology is relevant. It is known that for a large enough scale value $\alpha$, the persistence diagram $PD^{1}(VR_{\alpha}(\mathcal{G}))$ contains as many points (counting multiple points) as the genus of $\mathcal{G}$, the cardinality of a minimal system of cycles of $\mathcal{G}$~\cite{GGPSWWZ2017}. A main concern on the selection of landmarks to construct the witness or Dowker complexes is the loss of relevant topological information. Then, the next question arises:

\begin{problem}\label{pr:3} Let   $\mathcal{G}=(\mathcal{V},\mathcal{E},\omega)$ be weighted simple graph. Take the set of witnesses $W=\mathcal{V}$.
\begin{itemize}
    \item[(i)] Determine $\epsilon> 0$ and the way to select an $\epsilon$-net $L\subset \mathcal{V}$ such that $PD^{1}(Dow_{\alpha}(W,L))$ has cardinality  equal to the genus of $\mathcal{G}$, for a sufficiently large value $\alpha$ .
    \item[(ii)] Assuming that an $\epsilon$-net for $\mathcal{G}$ satisfying (i) is given. Find a bound for $d_{B}(PD^{1}(VR_{\alpha}(\mathcal{G})), PD^{1}(Dow_{\alpha}(W, L)))$, in terms of the lengths of cycles in $\mathcal{G}$ and $\epsilon$, where $d_{B}$ is the bottleneck distance in persistence diagrams.
\end{itemize}
 \end{problem} 

Finally, another important but largely unexplored topic
for PH obtained from the Dowker complexes is stability. 
Stability ensures that after applying a  ``small perturbation'' to  a metric graph $\mathcal{G}$ in order to obtain a new metric graph $\mathcal{G}'$, the corresponding persistence diagrams are close with respect to some distance, for instance the bottleneck or Wasserstein distance. By perturbation we refer to any transformation that may include edge deletion, addition, cleaving or contraction, or changing the weight function (which automatically modifies the distance $d_{\mathcal{G}}$). In the case of complexes depending on a selection of landmarks and witnesses a perturbation could be   changing the selection of these two sets. For the Dowker complex,  some stability results have been obtained, including the Dowker interleaving and the Dowker duality~\cite[see Section 4.2.3]{CDO2014}. 
In particular, Dowker duality can be interpreted as the property such that if $L,W\subset\mathcal{V}$ are the sets of landmarks and witnesses of $\mathcal{G}$, then $Dow_{\alpha}(W,L)$ and $Dow_{\alpha}(L,W)$ have the same homotopy type and, then, the two filtered simplicial complexes have the same persistent homology.

\begin{problem}\label{pr:4} Let $\mathcal{G}=(\mathcal{V},\mathcal{E},\omega)$ be a weighted graph and let $L,W\subset \mathcal{V}$ be sets of landmarks and witnesses for $\mathcal{G}$. 
\begin{itemize}
\item[(i)] Let $\mathcal{G}'=(\mathcal{V},\mathcal{E},\omega')$ be the metric graph with the same set of vertices and edges of $\mathcal{G}$, but another weight function $\omega'$. Find a bound for the bottleneck distance of the persistence diagrams corresponding to the Dowker complexes of $\mathcal{G}$ and $\mathcal{G}'$ in terms of the supremum distance between $\omega$ and $\omega'$. 
\item[(ii)] If two sets of landmarks $L$ and $L'$ for $\mathcal{G}$ are $\epsilon$-nets. Find conditions on the landmarks, such that $Dow_{\alpha}(W,L)$ and $Dow_{\alpha}(W,L')$ have the same homotopy type.
\end{itemize}
\end{problem}

Recent studies have addressed some  properties of the Dowker complexes in directed graphs. For instance, a complete characterization of the Dowker persistence diagrams for cycle networks is presented, and some stability properties have been proven for pair swaps by~\cite{CM2018} and for conjugate and  shift equivalent relations~\cite{Cote2023}. 
While both the theoretical properties and practical utility of the Dowker complexes on graphs yet remain largely unexplored, we believe that the Dowker complexes have a potential to address many bottlenecks that currently preclude broader applications of TDA in analysis of real-world graph datasets. 
\mySub{Acknowledgements}  This work was supported by the Office of Naval Research (ONR) award N00014-21-1-2530. Any opinions, findings, conclusions, or recommendations expressed in this paper are those of the authors and do not necessarily reflect the views of ONR.



\printbibliography[heading=subbibliography]
\end{refsection}

\begin{refsection}
\section{An Application of Probabilistic Combinatorics to Quantum Circuit Expressiveness}\label{sec:kay}
\begin{flushright}
{\it Bill Kay, Ryan Bennink}
\end{flushright}

\mySub{Introduction}
\label{sec:intro}
Quantum computing is a computing system in which certain physical properties are leveraged for computational advantages for some classes of problems. Recently, access to some quantum computing devices has grown. Many of the accessible devices are small and noisy, and the {\em variational approach} in which a traditional computer modifies parameters of some quantum circuit towards maximizing a function of its output has found success~\cite{mcclean2016theory,magann2021pulses,yuan2019theory}. 
The purpose of this document is to introduce a purely combinatorial problem, and to explain how it is related to a problem in quantum computing. We note that we are casting the analysis in~\cite{bennink} to a probabilistic setting regarding random binary vectors. Hence, progress in the combinatorial problem has an interpretation which impacts quantum computing.  We then present a sequence of examples, preliminary results, and open questions. 

Commutative quantum circuits, originally introduced as ``Instantaneous Quantum Polytime'' circuits \cite{Shepherd2009}, are a special class of quantum circuits that are relatively simple yet exhibit non-trivial quantum behavior.  For our purposes a commutative quantum circuit may be described as a circuit that first creates a uniform superposition of all $2^n$ computational states of $n$ qubits, then applies parity-dependent phases to selected subsets of qubits. Each subset is specified by a binary vector $a \in \{0,1\}^n$ and its associated phase is specified by a real angle $\theta$. 

\begin{definition}[Commutative Quantum Circuit]
\label{cqc}
A commutative quantum circuit specified by the binary vectors $a_1, \ldots,a_m \in \{0,1\}^n$ and associated parameters $\theta_1, \ldots, \theta_m \in \mathbb{R}$ creates a quantum state of the form
 \[
 \frac{1}{2^{n/2}} \left( \mathrm{e}^{\mathrm{i} \phi_{0 \ldots 0}}, \ldots, \mathrm{e}^{\mathrm{i} \phi_{1 \ldots 1}} \right) \in \mathbb{C}^{2^n}
 \]
\noindent where
\[
\phi_x = \sum_{i=1}^{m} \theta_i (-1)^{\langle a_i, x \rangle}
\]
and $\langle \cdot, \cdot \rangle$ denotes the inner product in GF(2)$^n$.
\end{definition}

In the context of variational quantum computing, the vectors $a_1,\ldots,a_m$ specify a fixed circuit structure and the angles $\theta_1, \ldots, \theta_m$ are the variational parameters.
A fundamental question concerns the \emph{expressiveness} of such: given the vectors $a_1, \ldots, a_m$, how much can the output state be varied by varying the angles  $\theta_1, \ldots, \theta_m$?  Intuitively, the expressiveness can be quantified by the expected (dis)similarity of two output states with randomly chosen parameter values: the more expressive the circuit, the more likely two random output states will be dissimilar.  As discussed in \cite{bennink}, for commutative quantum circuits this expectation value can be reduced to a combinatoric problem involving binary matrices described below. 
There it is shown that in the special case that the circuit is maximal, i.e.\ there is an independent phase for each nonempty subset of qubits, the expressiveness reduces to the problem of counting {\em abelian squares}.  Here we consider the more general problem of computing the expressiveness of commutative quantum circuits involving arbitrary collections of qubit subsets. In the next section we cast this problem in purely combinatorial language. We conclude with some examples, special cases, and open problems.



\mySub{A Combinatorial Problem}
\label{sec:acp}
In this section we will have several definitions where prescribed notations will aid in clarity. Here we fix several such notations for the section for ease of exposition. Fix $ n, t \in \mathbb{Z}^+$ and let $1 \leq m \leq 2^n - 1$. 

\begin{itemize}
    \item Let $\mathbf{x} = (x_1^T, \ldots, x_t^T)$ where $x_i \in \{0,1\}^n$ for $1 \leq i \leq t$.
    \item  Let $\mathbf{y} = (y_1^T, \ldots, y_t^T)$ where $y_i \in \{0,1\}^n$ for $1 \leq i \leq t$. 
    \item $\mathbf{A} = (a_1^T, \ldots, a_m^T)$ where $a_i \in \{0,1\}^n$ for $1 \leq i \leq m$. Moreover, the $\{a_i\}_{i=1}^m$ are non-zero and pairwise distinct.
\end{itemize}

For clarity,  $\mathbf{x}$ and $\mathbf{y}$ are $n \times t$ binary matrices with specified column labels, and $\mathbf{A}$ is an $n \times m$ binary matrix with {\em distinct, non-zero} columns with specified labels. For our specific application, we take the columns of $\mathbf{A}$ to be given in lexicographic order. We remark that each $a_i$ is a length $n$ binary vector and has a natural correspondence with subsets of indices of length $n$ vectors. 

We now need the following definition:
\begin{definition}[Expressiveness Indicator]
Given $\mathbf{A}, \mathbf{x}$ as above and $1 \leq i \leq m$, we define the $i$th {\em expressiveness indicator} $m_i^\mathbf{A} (\cdot)$ as:
\[
m_i^\mathbf{A} (\mathbf{x}) := \sum_{j = 1}^t \langle a_i, x_j\rangle_2
\]
where the inner product $\langle \cdot, \cdot\rangle_2$ is {\em mod } $2$ but the outer sum is integral. 
\end{definition}
When $\mathbf{A}$ is clear from context we simply write $m_i(\cdot)$. As these are not standard definitions, and there is a mix of modular and integral arithmetic,  we include a clarifying example:

\begin{example}
Let
\[
\mathbf{A} = 
\begin{pmatrix}
1 & 0 & 0 & 0&1&1&1\\
0 & 1 & 0 & 0&1&0&0\\
0 & 0 & 1 & 0&0&1&0\\
0 & 0 & 0 & 1&0&0&1

\end{pmatrix};
\mathbf{x} = 
\begin{pmatrix}
1 & 0 & 1 & 0\\
1 & 1 & 0 & 0\\
0 & 0 & 1 & 0\\
0 & 0 & 0 & 1\\

\end{pmatrix}
\]
 Then we have:
\begin{align*}
    m_5(\mathbf{x}) &= \langle a_5, x_1\rangle_2 + \langle a_5, x_2\rangle_2 +\langle a_5, x_3\rangle_2 +\langle a_5, x_4\rangle_2\\
    & = 0 + 1 + 1 + 0\\
    & = 2
\end{align*}
\end{example}

We are now ready to pose a general combinatorial question and relate it back to computing expressiveness of commutative  quantum circuits.

\begin{comp_problem}
Given $\mathbf{A}$,  let $\mathbf{x}$ and $\mathbf{y}$ be sampled uniformly at random. Compute:
\[
\mathbb{P}\left( (m_i(\mathbf{x}) = m_i(\mathbf{y})) \ \forall i \right) 
\]
\end{comp_problem}

While the Computational Problem is interesting in its own right, we would like to clarify how it is related to computing expressiveness of commutative quantum circuits. Recall from Definition~\ref{cqc} that the expressiveness of a commutative quantum circuit takes as input a collection of subsets of indices of a state vector. Further, we have remarked that for a length $n$ state vector, there is a natural correspondence between a subset of indices and binary length $n$ vectors (namely, the indicator vector of which indices are selected). Indeed, one can see in the analysis in~\cite{bennink} that if we encode the subsets of indices in the state vector we are interested in as $\mathbf{A}$, then the computation given in Computational Problem is the critical computation for finding the desired expressiveness. In particular, it is shown that expressiveness of {\em maximal} commutative quantum circuits is equivalent to taking $\mathbf{A}$ to be all $2^n -1$ non-zero length $n$ binary vectors, and that the probability in Computational Problem is precisely resolved by counting abelian squares.

\mySub{Examples, Preliminary Results, and Open Problems}

Notice that for each choice of $\mathbf{A}$, there is a different instance of Computational Problem. For some fixed choices of $\mathbf{A}$, there is a natural interpretation of Computational Problem. 

\begin{example}
Let $\mathbf{A}$ be the collection of all weight 1 length $n$ vectors (i.e., the $n \times n$ identity matrix). Let $\mathbf{x}$ and $\mathbf{y}$ be $n \times t$ binary matrices. Then 
\[
 (m_i(\mathbf{x}) = m_i(\mathbf{y})) \ \forall i 
\]
if and only if $\mathbf{x}$ and $\mathbf{y}$ have the same row sums. 
\end{example}

This is clear, as $\langle a_j, x_i\rangle_2 = 1$ if and only if $x_i$ is $1$ in coordinate $j$. Using this combinatorial interpretation, we can resolve Computational Problem when $\mathbf{A}$ is the $n \times n$ identity matrix  by the following (equivalent) lemma:

\begin{lemma}
\label{l1}
Let $\mathbf{x}$ and $\mathbf{y}$ be $n \times t$ binary matrices sampled uniformly at random. Then the probability that $\mathbf{x}$ and $\mathbf{y}$ have all the same row sums is 
\[
\binom{2t}{t}^n\frac{1}{2^{2tn}}.
\]
\end{lemma}

\begin{proof}
We will argue for a single pair of uniform random binary vectors of length $t$. Since we are computing row sums for $\mathbf{x}$ and $\mathbf{y}$, we argue first  for a fixed row and conclude by taking a product of independent probabilities. Summing over possible row sums, the probability that a pair of length $t$ binary vectors have the same row sum is given by:
\begin{align*}
\sum_{k=0}^t \binom{t}{k}^2 (1/2)^{2k}(1/2)^{2(t-k)} &=\frac{1}{2^{2t}}\sum_{k=0}^t \binom{t}{k}^2 \\
& = \frac{1}{2^{2t}}\sum_{k=0}^t\binom{t}{k}\binom{t}{t-k}\\
& = \frac{1}{2^{2t}} \binom{2t}{t}
\end{align*}
 where the last equality can be seen by counting the number of length $2t$ binary vectors with precisely $t$ ones. The right hand side counts this quantity directly, while the left hand side counts the number of ways for $k$ ones in the first half with $t-k$ ones in the second half. Taking the product over all $n$ (independent)  rows of $\mathbf{x}$ and $\mathbf{y}$ gives the desired probability. 
\end{proof}

Each column of $\mathbf{A}$ provides an additional constraint on possible choices of $\mathbf{x}$ and $\mathbf{y}$ in Computational Problem. In fact, the collections of $(\mathbf{x},\mathbf{y})$ pairs which satisfy the computation in Computational Problem are {\em monotone} in the columns of $\mathbf{A}$, in that if the set of columns $\mathbf{A}_2$ contains the set of columns of  $\mathbf{A}_1$, then 
\[
 (m_i^{\mathbf{A}_2}(\mathbf{x}) = m_i^{\mathbf{A}_2}(\mathbf{y})) \ \forall i \Rightarrow (m_i^{\mathbf{A}_1}(\mathbf{x}) = m_i^{\mathbf{A}_1}(\mathbf{y})) \ \forall i.
\]

For example, if  $\mathbf{A}$ has among its  columns all weight $1$ length $n$ vectors, if $\mathbf{x}$ and $\mathbf{y}$ satisfy the condition in Computational Problem then $\mathbf{x}$ and $\mathbf{y}$ necessarily have the same row sums. We provide combinatorial interpretation of Computational Problem for another choice of $\mathbf{A}$, and then exploit the monotonicity property to state and prove Computational Problem for another choice of $\mathbf{A}$. 

\begin{example}
\label{e2}
Let $a_{i,j}$ be the length $n$ binary vector with a $1$ in positions $i$ and $j$, $0$ otherwise. Let $x$ be any length $n$ binary vector. $\langle a_{i,j}, x \rangle_2$ is $1$ if the $i$th and $j$th entries of $x$ disagree, $0$ otherwise. We call this an {\em $i,j$-mismatch} in $x$. Hence, if $\mathbf{A}$ consists of all length $n$ binary vectors of weight at most $2$, and  $\mathbf{x}$ and $\mathbf{y}$ are $n \times t$ binary matrices, then 
\[
 (m_i(\mathbf{x}) = m_i(\mathbf{y})) \ \forall i 
\]
if and only if $\mathbf{x}$ and $\mathbf{y}$ have the same row sums, and for each pair of rows $i$ and $j$  the number of columns in which $\mathbf{x}$ and $\mathbf{y}$ have an $i,j$-mismatch are the same.
\end{example}

Computational Problem for the family presented in Example~\ref{e2} remains open. However, we present the following:

\begin{lemma}
\label{l2}
Let $\mathbf{A}$ be the collection of binary vectors weight $1$ and weight $2$ with a $1$ in the first coordinate. Let $\mathbf{x}$ and $\mathbf{y}$ be $n \times t$ binary matrices sampled uniformly at random. Then:
\[
\mathbb{P}\left( (m_i(\mathbf{x}) = m_i(\mathbf{y})) \ \forall i \right) 
\]

is given by:

\[
\frac{1}{2^{2tn}} \sum_{k = 0}^t \left(\binom{2k}{k} \binom{2(t-k)}{t-k}\right)^{n-1}\binom{t}{k}^2
\]
\end{lemma}

\begin{proof}
 For $\mathbf{x}$ and $\mathbf{y}$ sampled uniformly at random, we want to compute the probability  that $\mathbf{x}$ and $\mathbf{y}$ have the same row sums, and for each $ 1 < j \leq n$,  the number columns which have a $1,j$-mismatch is the same. We claim that 
 \[
\sum_{k=0}^t \binom{t}{k}^2\left(\sum_{j = 0}^k \sum_{\ell = 0}^{t-k} \binom{k}{j}^2 \binom{t-k}{\ell}^2\right )^{n-1} 
 \]
counts the number of $(\mathbf{x}, \mathbf{y})$ pairs for which 
\[
 (m_i(\mathbf{x}) = m_i(\mathbf{y})) \ \forall i  
\]
 holds. We now count such pairs. First, fix a first row for each of $\mathbf{x}$ and $\mathbf{y}$. Since the row sums must be the same, these first rows have some number of $1$s, say $k$. For each row $1< r\leq n$, there are some number of $1$s in $\mathbf{x}$ and $\mathbf{y}$. Since the row sums are the same for $\mathbf{x}$ and $\mathbf{y}$, this number is the same. We note that $\mathbf{x}$ and $\mathbf{y}$ need to have the same number of $1,r$ matches as well as mismatches, and this can be counted by enumerating the number of ways to place $1$s in columns with a $1$ in the first entry, and then columns with $0$ in the first entry. There are $j + \ell$ 1s in row $r$, with $j$ in columns with $1$ in the first entry and $\ell$ in columns with $0$ in the first entry. Hence, for each fixed choice of $j$ and $\ell$ there are
\[
 \binom{k}{j}^2 \binom{t-k}{\ell}^2
\]
ways to place these $1$s. Summing over choices of $j$ and $\ell$ we see that for each fixed pair of first rows for $\mathbf{x}$ and $\mathbf{y}$ there are
\[
\sum_{j=0}^k \sum_{\ell = 0}^{t-k} \binom{k}{j}^2 \binom{t-k}{\ell}^2
\]
choices to place the $1$s in row $r$. These choices are independent for rows $1 < r \leq n$, and so there are 
\[
\left(\sum_{j=0}^k \sum_{\ell = 0}^{t-k} \binom{k}{j}^2 \binom{t-k}{\ell}^2\right)^{n-1}
\]
choices for $\mathbf{x}$ and $\mathbf{y}$ for a fixed first row. For each $k$, there are $\binom{t}{k}$ ways to place $1$s in the first row of $\mathbf{x}$ and $\mathbf{y}$. Summing over each $k$, the total number of desired $(\mathbf{x},\mathbf{y})$ pairs is
\[
\sum_{k=0}^t \binom{t}{k}^2\left(\sum_{j=0}^k \sum_{\ell = 0}^{t-k} \binom{k}{j}^2 \binom{t-k}{\ell}^2\right)^{n-1}
\]

Applying analysis similar to that in the proof of Lemma~\ref{l1} shows that this is 
\[
\sum_{k = 0}^t \left(\binom{2k}{k} \binom{2(t-k)}{t-k}\right)^{n-1}\binom{t}{k}^2
\]
There are $2^{2tn}$ total possible pairs $(\mathbf{x},\mathbf{y})$, and dividing through yields the claim in the statement of Lemma~\ref{l2}. 


\end{proof}

The methods in the proofs of Lemmas~\ref{l1} and~\ref{l2} depend on row-wise independence. However, these methods do not directly extend to $\mathbf{A}$ for, say, all  binary vectors of weight at most $2$, as for $1 \leq i < j< k \leq m$, the configurations of $i,j$-mismatches and the configurations of $j,k$-mismatches affect the possible configurations of $i,k$- mismatches. More broadly, once the family $\mathbf{A}$ becomes rich enough, the step where we appeal to independence is not applicable and more sophisticated analysis is required. Hence, we have the following (open) questions:

\begin{problem}
\label{q1}

Fix $m$ and $k \leq m$. Can we answer Computational Problem when $\mathbf{A}$ consists of all vectors of weight at most $k$? 

\end{problem}

Lemma~\ref{l1} answers Open Problem~\ref{q1} when $k =1$, and the maximal case in~\cite{bennink} answers Open Problem~\ref{q1} when $k=m$. These boundary cases are the only values of $k$ for which Open Problem~\ref{q1} has been answered.    

\begin{problem}
\label{q2}
For which choices of $\mathbf{A}$ can we answer Computational Problem exactly? For which choices of $\mathbf{A}$ can we answer Computational Problem efficiently?
\end{problem}

To rephrase the second part of Open Problem~\ref{q2}, we ask:

\begin{problem}
\label{q3}
Is Computational Problem NP-Complete for arbitrary $\mathbf{A}$?
\end{problem}

Finally, we remark that exact answers to Computational Problem provide answers to expressiveness of commutative quantum circuits, approximate answers to Computational Problem provide approximations to expressiveness of commutative quantum circuits. Hence we ask:

\begin{problem}
\label{q4}
For which choices of $\mathbf{A}$ can we approximate answers to Computational Problem?
\end{problem}
\mySub{\bf Acknowledgements}
This work was performed at Oak Ridge National Laboratory, operated by
UT-Battelle, LLC for the US Department
of Energy (DOE) under contract DE-AC05-00OR22725. Support for the work came from the DOE Advanced
Scientific Computing Research (ASCR) Accelerated Research in Quantum
Computing Program under field work proposal ERKJ354.

\printbibliography[heading=subbibliography]
\end{refsection}

\begin{refsection}
\section{The Computational Complexity of the 4-uniform Hypergraph Minimum $s$-$t$ Cut Problem}\label{sec:veldt}
\begin{flushright}
	{\it Nate Veldt}
\end{flushright}

\begin{figure}[h]
	\centering
	\begin{minipage}{.35\textwidth}
		\centering
		\includegraphics[width=\linewidth]{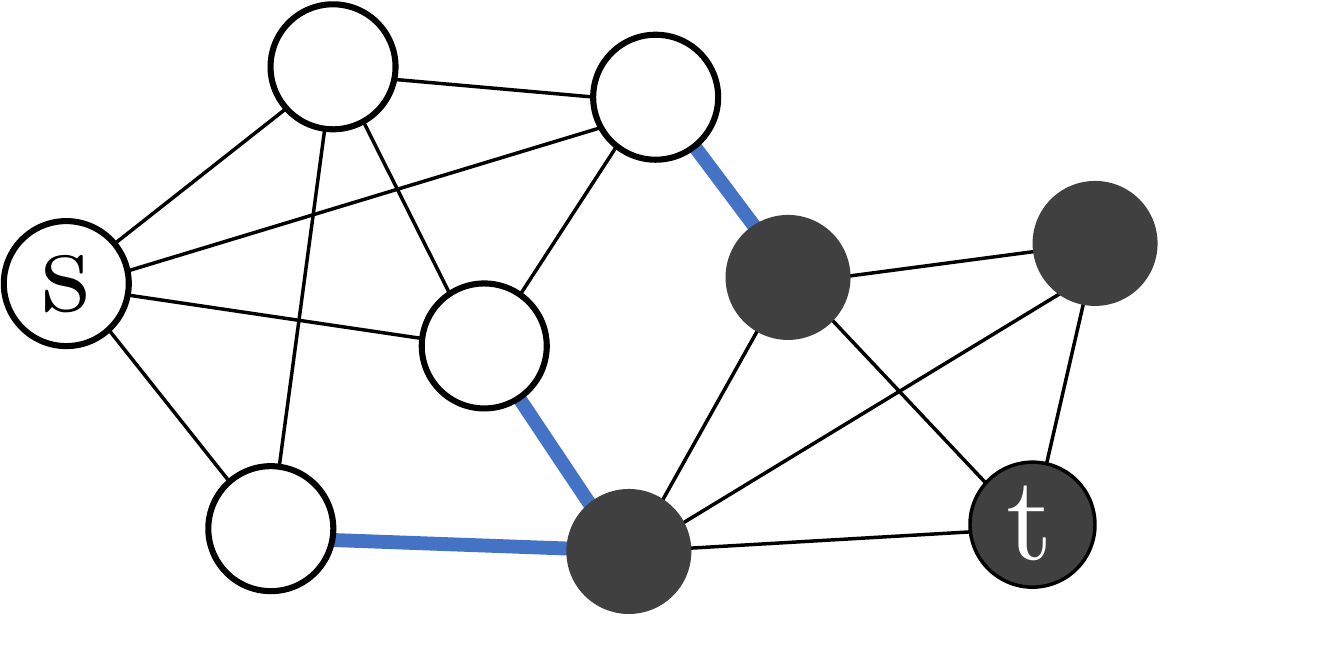}
		\caption{Graph $s$-$t$ cut.}
		\label{fig:graph}
	\end{minipage}%
	\begin{minipage}{0.65\textwidth}
		\textbf{Two formulations of the graph $s$-$t$ cut problem}
		\begin{enumerate}
			\item Find a minimum size set of edges whose deletion separates $s$ from $t$ (e.g., three blue edges in Figure~\ref{fig:graph}).
			\item Partition $V$ into $S \subseteq V$ and $\bar{S} = V\backslash S$ (e.g., white and black nodes in Figure~\ref{fig:graph}), with $s \in S$ and $t \in \bar{S}$, in order to minimize the number of \emph{cut} edges (i.e., edges crossing the partition).
		\end{enumerate}
	\end{minipage}
\end{figure}

\mySub{Introduction} 
Finding a minimum $s$-$t$ cut in a graph is one of the most fundamental and well-known combinatorial optimization problems, studied and applied widely in mathematics, computer science, operations research, and data science. The input to the problem is a graph $G = (V,E)$ with a distinguished source node $s$ and sink node $t$. The goal is to delete the minimum number of edges to destroy all paths from $s$ to $t$, or equivalently partition the nodes in a way that separates $s$ and $t$ while minimizing a cut penalty (see Figure~\ref{fig:graph}). A polynomial time solution for this problem has been known since the 1950s~\cite{ford1956maximal}, and many other algorithms are now taught regularly as a standard part of an undergraduate algorithms course. 

\begin{figure}[h]
	\centering
	\includegraphics[width=.9\linewidth]{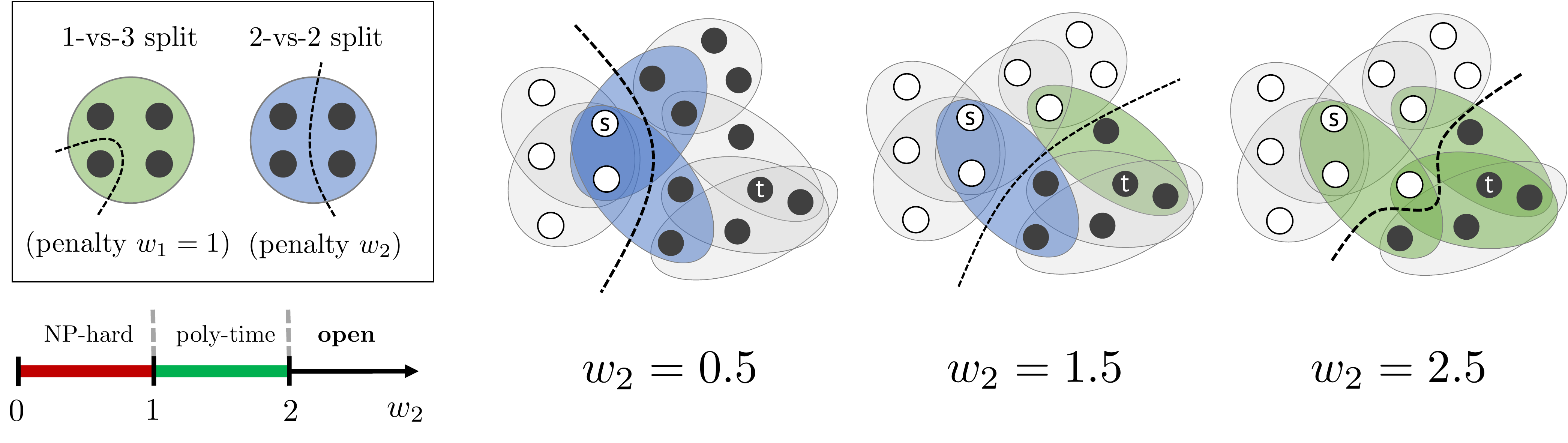}
	\caption{Optimal $s$-$t$ cuts in a 4-uniform hypergraph, depending on the penalty $w_2$ assigned to $2$-vs-$2$ hyperedge splits (i.e., two nodes on each side of a bipartition). All $1$-vs-$3$ hyperedge splits have a penalty $w_1 = 1$. Previous results have shown the problem is NP-hard for $w_2 \in [0,1)$ and polynomial-time solvable when $w_2 \in [1,2]$. The complexity of the problem remains open for $w_2 > 2$.}
	\label{fig:4unif}
\end{figure}

The $s$-$t$ cut problem can also be defined when the input is a hypergraph $\mathcal{H} = (V,E)$, in which case each edge $e \in E$ can have two \emph{or more} nodes. If we generalize one way of formulating the graph $s$-$t$ cut problem, the goal is simply to delete a minimum sized set of hyperedges to separate a source node $s$ from a sink node $t$. A polynomial time solution for this problem was given half a century ago by Lawler~\cite{lawler1973}, which works by reducing the hypergraph to an $s$-$t$ cut problem in a \emph{directed} graph with an augmented node set. However, generalizing the second formulation of the graph $s$-$t$ cut problem (partition nodes to separate $s$ and $t$ while minimizing a {cut} penalty) leads to alternative problems that are not equivalent to this edge-deletion problem. When a hyperedge has more than two nodes, there is more than one way to separate those nodes across two node sets $S$ and $\bar{S}$. In a 4-uniform hypergraph (see Figure~\ref{fig:4unif}), we can treat $1$-vs-$3$ hyperedge splits differently from $2$-vs-$2$ splits. Recently there has been a growing interest in solving hypergraph cut problems under these generalized types of cut penalties~\cite{veldt2020hyperlocal,veldt2021approximate,veldt2022hypergraph,zhu2022hypergraph,panli_submodular}, since different ways of splitting up the nodes of a hyperedge may be more or less desirable depending on the application. In many machine learning and data mining applications, for example, hyperedges represent evidence that a set of data objects are related and should be associated with the same label or cluster. In these cases, it is more desirable to split hyperedges in such a way that most (even if not all) nodes are on the same side of a cut. At the same time, using different cut penalties leads to drastic differences in the underlying computational complexity of the problem, even when hyperedges are very small. We specifically consider an open question on the computational complexity of a generalized hypergraph $s$-$t$ cut problem in 4-uniform hypergraphs.

\mySub{The cardinality-based hypergraph $s$-$t$ cut problem}
Let $\mathcal{H} = (V,\mathcal{E})$ be a hypergraph with a source node $s$ and sink node $t$. Given a bipartition $\{S, \bar{S}\}$ of the nodes, hyperedges can be classified based on the number (i.e., cardinality) of nodes on the small side of the cut. Edges with $i$ nodes on the small side of the cut are denoted by $\partial_i(S) = \{ e \in \mathcal{E} \colon \min \{|e \cap S|, |e \cap \bar{S} |\} = i\}$.
We take the minimum between $|e \cap S|$ and $|e \cap \bar{S}|$ to ensure the resulting hypergraph cut function is symmetric, generalizing the fact that a graph cut function is symmetric.
Let $w_i \geq 0 $ be the penalty assigned to each hyperedge in $\partial_i(S)$. The value $w_0$ is set to 0 to ensure that a hyperedge only has a penalty if it is actually {cut}. 
The {cardinality-based} hypergraph $s$-$t$ cut problem is given by
\begin{equation}
	\label{eq:cbstcut}
	\min_{S \subseteq V} \;\; \sum_i w_i |\partial_i(S) |, \hspace{.5cm} \text{ subject to $s \in S$ and $t \in \bar{S}$}.
\end{equation}
When $w_i = 1$ for every $i > 0$, this is equivalent to finding a minimum number of hyperedges to remove to separate $s$ from $t$, and can be solved in polynomial time using the reduction of Lawler~\cite{lawler1973}. If $\mathcal{H}$ is 3-uniform, any way of cutting a hyperedge places exactly one node in  $S$  or exactly one node in $\bar{S}$. Since all hyperedge cut penalties are the same, this is solved by Lawler's algorithm.


\mySub{Results for 4-uniform hypergraphs}
When $\mathcal{H}$ is 4-uniform, there is a distinction between 2-vs-2 splits (penalty of $w_2$) and 1-vs-3 splits (penalty of $w_1$). As long as $w_1 > 0$, parameters can be scaled without loss of generality so that $w_1 = 1$. For certain choices of $w_2$, Veldt, Benson, and Kleinberg~\cite{veldt2022hypergraph} showed how to reduce the hypergraph problem to a graph $s$-$t$ cut problem, by replacing each hyperedge with a \emph{gadget} involving auxiliary vertices and directed edges. The aim is to design gadgets so that cut penalties in the reduced directed graph match cut penalties in the original hypergraph. In the 4-uniform case, it turns out this is possible if \emph{and only if} $1 \leq w_2 \leq 2$. If $0 \leq w_2 < 1$, the problem can be shown to be NP-hard via reduction from the optimization version of maximum cut, one of Karp's 21 NP-complete problems~\cite{karp2010reducibility}. If $G = (V,E)$ is an unweighted graph representing an instance of maximum cut, this reduction constructs a hypergraph with the same node set plus two additional nodes $s$ and $t$. For each edge $(u,v) \in E$, a hyperedge $\{u,v,s,t\}$ is introduced, which must be cut in the resulting hypergraph $s$-$t$ cut problem since it contains both $s$ and $t$. If $w_2 < w_1 = 1$, then a 2-vs-2 split of the hyperedge is cheaper than a 1-vs-3 split. Therefore, the minimum $s$-$t$ cut is the cut maximizing the number of 2-vs-2 hyperedge splits, which is equivalent to finding a bipartition of nodes in $G$ that maximizes the number of cut edges.

\mySub{Open questions and motivation}
Given the above notation and terminology, our open question can be stated formally as follows.
\begin{problem}
	\label{op:4hyper}
	What is the complexity of the cardinality-based 4-uniform hypergraph $s$-$t$ cut problem when $w_2 > 2$ and $w_1 = 1$?
\end{problem}
Nothing is known about the computational complexity of the problem except that the graph reduction strategy no longer applies. A reasonable first step is to address the special case obtained by setting $w_1 = 1$ and taking the limit as $w_2 \rightarrow \infty$. This limit converges to the \textsc{No-Even-Split} 4-uniform $s$-$t$ cut problem: minimize the number of $1$-vs-$3$ hyperedge splits when separating $s$ and $t$, without making any $2$-vs-$2$ hyperedge splits. Placing node $s$ in a cluster by itself always provides one way to separate $s$ from $t$ without having even splits, but the complexity of finding the \emph{minimum} number of 1-vs-3 splits is unknown.
\begin{problem}
	What is the complexity of the \textsc{No-Even-Split} 4-uniform $s$-$t$ cut problem?
\end{problem}
This special case of Open Problem~\ref{op:4hyper} is particularly interesting given its close relationship to an $s$-$t$ cut problem with a very simple solution. More specifically, this problem at first appears equivalent to setting $w_2 = 1$ and taking a limit as $w_1 \rightarrow 0$, but there is a subtle and interesting difference. Setting $w_2 = 1$ and $w_1 = 0$ results in a degenerate problem that is easily optimized, since separating $s$ from the rest of the nodes has a cut penalty of 0. However, if $w_1 = 1$ and $w_2 \rightarrow \infty$, it may never be optimal to separate $s$ (or $t$) from all other nodes. Despite its close relationship to a degenerate and easily solved problem, it is not clear if \textsc{No-Even-Split} 4-uniform $s$-$t$ cut is NP-hard, or if there is a strategy leading to a polynomial time solution.

A solution to the \textsc{No-Even-Split} problem would hopefully shed light on other polynomial time algorithms or hardness results for the parameter region $w_2 \in (2,\infty)$. In hypergraphs where the maximum hyperedge size is greater than 4, there is an even larger gap between known NP-hardness results and polynomial time solutions. When hyperedges are of arbitrary size and cut penalties satisfy $w_i = f(i)$ for some increasing concave function $f$, the hypergraph problem can be reduced to a directed graph $s$-$t$ cut problem~\cite{veldt2022hypergraph}. There are also some parameter regimes where the problem is known to be NP-hard, but there are many parameter settings where neither hardness results nor polynomial time solutions are known. New techniques for the 4-uniform case would hopefully help close these gaps as well. 

In addition to their theoretical value, answers to these questions have the potential to advance the state of the art in practical hypergraph algorithms for downstream applications. Generalized hypergraph cut problems are already being used in data mining and machine learning applications such as node classification~\cite{panli_submodular,veldt2020hyperlocal,zhu2022hypergraph} and image segmentation~\cite{veldt2021approximate}, but existing methods focus only on the regime where cut penalties are submodular. Cut penalties are also cardinality-based in most applications~\cite{veldt2022hypergraph}, in which case this submodularity property exactly corresponds to cut penalties defined by a concave function $f$. If, for example, penalties were instead chosen so that $w_i = g(i)$ for an increasing and \emph{strictly convex} function $g$, this would provide even more incentive to assign nodes from the same hyperedge to the same cluster or node label. However, the computational complexity of finding minimum cuts in this setting is unknown. The $4$-uniform hypergraph $s$-$t$ cut problem with $w_2 > 2$ is the simplest example of this type of convex-penalty hypergraph cut problem. An answer to Open Problem~\ref{op:4hyper} would provide a needed first step in understanding what is algorithmically possible for this and many other types of generalized hyperedge cut penalties.

\printbibliography[heading=subbibliography]
\end{refsection}

\begin{refsection}
\section{The Edge and Vertex Elimination Problems in Directed Acyclic Graphs}\label{sec:naumann1}

\begin{flushright}
{\it Uwe Naumann}
\end{flushright}
\date{}

\mySub{Introduction}

Impressive progress in the development of computer hard- and software has been 
made over the past decades. Consequently, numerical simulation has become one 
of the pillars of science and engineering. Practically relevant real-world 
phenomena are modelled mathematically. Their numerical evaluation yields 
multivariate vector functions 
$F: \R^n \rightarrow \R^m$
implemented as often 
highly complex computer programs. 

Algorithmic Differentiation \cite{Griewank2008EDP,Naumann2012TAo} of numerical 
programs plays a central role in numerous
areas of computational science and engineering including error estimation, uncertainty quantification, parameter sensitivity analysis, model calibration and 
optimization. 
The resulting derivative programs
can be used, for example, to compute the Jacobian (matrix) 
of (a differentiable) $F.$ 

The data flow in a numerical program induces a directed acyclic graph (dag)
$G=(V,E)$ with integer vertices representing all input, intermediate and output values $v_i,$ $i=1,\ldots,|V|$ of the program and directed edges 
$E \subseteq V \times V$ modelling data dependence.
Association of local partial derivatives 
$d_{i,j}=\frac{\partial v_j}{\partial v_i}$ 
with all edges yields the chain rule of differentiation 
as
\begin{equation} \label{eqn:Baur}
\frac{d v_t}{d v_s} = \sum_{(s,*,t)} \prod_{(i,j) \in (s,*,t)} d_{i,j} \; ,
\end{equation}
where summation is over all paths $(s,*,t)$ connecting a vertex $s$ with 
a vertex $t$ \cite{Baur1983TCo}. (Partial) Derivatives are thus defined for 
arbitrary pairs of values represented by both vertices. 

\mySub{Formal Statement of the Problem}

\begin{definition}[Edge Elimination] \label{ee}
	Let $G=(V,E)$ and $(i,j) \in E.$
	{\em Front-elimination} of $(i,j)$ yields
	$G - (i,j) \equiv G^f=(V^f,E^f)$ such that
	$$
	V^f = \begin{cases}
		V & \text{if}~|P_i|>1 \\
		V \setminus j & \text{otherwise} \\
	\end{cases} \quad \text{and} \quad 
	E^f= E \, \cup \, \{(i,k) : k \in S_j\} 
	$$
at the cost of $|S_j|.$

	{\em Back-elimination} of $(i,j)$ yields
	$G - (j,i) \equiv G^b=(V^b,E^b)$ such that
	$$
	V^b = \begin{cases}
		V & \text{if}~|S_i|>1 \\
		V \setminus i & \text{otherwise} \\
	\end{cases} \quad \text{and} \quad 
	E^b= E \, \cup \, \{(k,j) : k \in P_i \}
	$$
at the cost of $|P_i|.$
\end{definition}
{\em Front-[back-]eliminatable} edges are required to have non-empty 
successor [predecessor] sets. By the chain rule of differentiation, all 
{\em complete} edge elimination sequences transform $G$ into a bipartite 
dag representing the Jacobian matrix of the underlying numerical program.
The cost of an edge elimination sequence is the sum of the costs of the 
individual edge eliminations.
We aim to minimize the cost over all complete edge elimination sequences.
\begin{problem}[{\sc Edge Elimination}]
Given a directed acyclic graph $G$ and a positive integer $k \geq 0,$ is there a complete edge
elimination sequence with cost less than or equal to $k?$
\end{problem}
The elimination of a vertex is equivalent to 
front-elimination of its in-edges. 
Similarly, it is equivalent to back-elimination of its out-edges.
\begin{definition}[Vertex Elimination] \label{ve}
	Let $G=(V,E)$ and $j \in V.$
	{\em Elimination} of $j$ yields
	$G - j \equiv G'=(V',E')$ such that
	$$
	V' = V \setminus j 
	\quad \text{and} \quad 
	E'= E \, \cup \, \{ (i,k) : i \in P_j~\text{and}~k \in S_j \}
	$$
at the cost of $|P_j| \cdot |S_j|.$
\end{definition}
{\em Eliminatable} vertices are required to have non-empty 
predecessor and successor sets. By the chain rule of differentiation, all 
{\em complete} vertex elimination sequences transform $G$ into a bipartite 
dag representing the Jacobian matrix of the underlying numerical program.
The cost of a vertex elimination sequence is the sum of the costs of the 
individual vertex eliminations.
We aim to minimize the cost over all complete vertex elimination sequences.
\begin{problem}[{\sc Vertex Elimination}]
Given a directed acyclic graph $G$ and a positive integer $k \geq 0,$ is there a complete vertex
elimination sequence with cost less than or equal to $k?$
\end{problem}

\mySub{What Is [Not] Known}

Both edge and vertex elimination terminate. Structural and numerical 
correctness follows by the chain rule \cite{Naumann2004Oao}.

\begin{theorem}
{\sc Jacobian Accumulation} is NP-complete.
\end{theorem}
The proof, see \cite{Naumann2008OJa}, uses reduction from {\sc Ensemble 
Computation} \cite{Garey1979CaI}. 
It exploits potential algebraic dependences (in particular, 
equality) among local partial derivatives (labels on edges in the dag).
Under the same assumptions both {\sc Edge Elimination} and 
{\sc Vertex Elimination} turn out to be NP-complete. The computational
complexity of the purely structural formulations in Defs.~\ref{ee} and \ref{ve}
is unknown.

Exiting heuristics for computing good vertex and edge elimination sequences are 
based on the structural formulations, for example \cite{Chen2012AIP,Forth2004JCG,Griewank2003AJa,Pryce2008FAD}.

Incomplete elimination sequences are required in the case of {\em scarcity}
\cite{Griewank2005AaE}. The computational complexity of the 
associated combinatorial {\sc Minimum Edge Count} 
problem asking for a (matrix-free) representation of the Jacobian as a dag 
with a minimal number of edges is unknown \cite{Mosenkis2012Oop}. Ultimately,
combinations {\sc Edge Elimination} and {\sc Vertex Elimination} with 
{\sc Minimum Edge Count} should be considered.
\printbibliography[heading=subbibliography]
\end{refsection}

\begin{refsection}
\section{Data Flow Reversal Problems}\label{sec:naumann2}

\begin{flushright}
{\it Uwe Naumann}
\end{flushright}

\mySub{Introduction}

Impressive progress in the development of computer hard- and software has been 
made over the past decades. Consequently, numerical simulation has become one 
of the pillars of science and engineering. Practically relevant real-world 
phenomena are modelled mathematically. Their numerical evaluation yields 
multivariate vector functions 
$F: \R^n \rightarrow \R^m :$ $\Y=F(\X)$
implemented as often 
highly complex computer programs. 

Evaluation of a numerical program at a given input induces a directed 
acyclic graph (DAG) $G=(V,E).$ 
Its integer vertices $V=(X,Z,Y)$ consist of $n \geq 1$ sources 
$X=\{1,\ldots,n\}$ representing $\X \in \R^n,$
$p \geq 0$ intermediate vertices $Z=\{n+1,\ldots,n+p\}$ and $m \geq 0$ 
sinks $Y=\{n+p+1,\ldots,n+p+m\}$ representing $\Y \in \R^m.$ Intermediate 
vertices and sinks model
calls to elemental functions
$$
\V_j=F_j(\V_i)_{(i,j) \in E} \quad \text{for}~j=n+1,\ldots,n+p+m \; .
$$
Directed edges $E \subseteq V \times V$ are due to corresponding 
data dependence within the 
program. A topological order of the vertices is implied. 
DAGs of numerical simulations become very large.
They are nonpersistent in the sense that they cannot be stored in memory.

Algorithmic Differentiation (AD) \cite{Griewank2008EDP,Naumann2012TAo} of numerical 
programs in adjoint mode plays a central role in numerous
areas of computational science and engineering including error estimation, uncertainty quantification, parameter sensitivity analysis, model calibration and 
optimization. It yields efficient (${\mathcal O}(m)$) gradients ($m=1$) with 
machine accuracy which would otherwise have to be approximated numerically 
with a cost of ${\mathcal O}(n)$ where often $n \gg 1.$ Refer to 
\url{www.autodiff.org} for links to AD projects and research groups as well as
for an extensive bibliography on the subject.

The {\em algorithmic adjoint} evaluates 
$\bar{\X}=\bar{\X}+F'(\X)^T \cdot \bar{\Y}$ as
$$
\bar{\V}_i=\bar{\V}_i+F'_j(\V_i)_{(i,j) \in E} \cdot \bar{\V}_j \quad \text{for}~j=n+p+m,\ldots,n+1 
$$
by accessing the nonpersistent $\V_i$ in reverse order. 
Equivalently, the DAG $G$ is reversed by accessing its vertices in 
reverse order.
Values lost due to overwriting of program
variables in the given implementation of $F$ need to be recovered. Naively,
all overwritten {\em required} \cite{Dauvergne2006TDF,Hascoet2005TBR} 
values could be pushed onto a stack. The persistent memory required for data 
flow reversal becomes maximal (${\mathcal O}(n+p)$). Alternatively, 
all overwritten required values could be recomputed 
from the persistent inputs, which yields minimal persistent memory requirement 
(${\mathcal O}(n)$). 
Reversal of the flow of data in nonpersistent memory results in quadratic 
computational cost (${\mathcal O}(p^2)$). Each intermediate value needs to be
recomputed from $\X.$ Without loss of generality, unit cost can be
assumed for all $F_j.$
The combinatorial {\sc DAG Reversal} 
problem aims to minimize the computational cost of data flow reversal in 
limited ($\ll {\mathcal O}(p^2)$) persistent memory.

Applications of data flow reversal beyond algorithmic adjoints include 
reverse debugging of large-scale numerical simulations \cite{Engblom2012ARO} 
and checkpointing for resilience of computer systems \cite{osti_1364654}.

\mySub{Formal Statement of the Problem}

\begin{problem}[{\sc DAG Reversal}] \label{dagr}
Given a DAG $G$ and two integers $M \geq n$ and $O \geq 0$ 
can $G$ be reversed within persistent memory of size less than or equal to
$M$ with computational cost less than or equal to $O?$
\end{problem}

\mySub{What Is [Not] Known}

\begin{theorem}
{\sc DAG Reversal} is NP-complete.
\end{theorem}
The proof can be found in \cite{Naumann2008DRi}.
The related {\sc Call Tree Reversal} problem is also known to be NP-complete
\cite{Naumann2008CTR}.

The state of the art in methods for data flow reversal is mostly driven by 
solutions for special cases due to given simulation scenarios 
\cite{Charpentier2001CSo,Symes2007Rtm,Wang2009MRD}. 
Optimal checkpointing methods for certain types of evolutions (DAG is a simple chain) were developed \cite{Griewank1992ALG,Walther2000Rev}. Generalizations
include mixed-integer programming for {\sc Call Tree Reversal} \cite{Lotz2016MIP} and divide-and-conquer approaches \cite{Siskind2018Dac}. Refer to
\cite{Hascoet2009RSf} for additional information.

The combinatorial data flow reversal problems are expected to benefit
from further formal theoretical analysis. In particular, there are no 
approximation methods available.
\printbibliography[heading=subbibliography]
\end{refsection}

\begin{refsection}




\section{Price of Asynchrony}\label{sec:ratfish}
\begin{flushright}
{\it Carlos Ortiz Marrero, Stephen J.\ Young}
\end{flushright}
\mySub{Introduction}
Finding the solution $x^*$ of a large-scale system of linear equations $Ax = b$ is one of the most common computation primitives in scientific computation, forming the backbone of numerical approaches  to a variety of applications including computational fluid dynamics, optimization, and atmospheric modeling.  While there are many different approaches to solve such systems (for example, Gaussian elimination, Krylov subspace methods, or conjugate gradient methods), iterative approaches are popular because of their straight-forward nature and ease of implementation.  Broadly speaking, iterative methods proceed by cleverly selecting an auxiliary matrix $M$ and using $M$ to generate a sequence of approximate solutions $x^{(0)}, x^{(1)}, x^{(2)}, \ldots$ with the recurrence $Mx^{(k+1)} = (M-A)x^{(k)} + b$.  Thus if systems of the form $My = c$ can be easily solved, then the approximation can be refined for essentially the cost of one matrix-vector multiplication.  For example, the (damped) Jacobi iteration is when $M$ is chosen to be the (scaled) diagonal of $A$ and Gauss-Seidel iteration is when $M$ chosen to be the lower-triangular portion of $M$.

Typically, the analysis of these methods work with the sequence of residuals $r^{(k)} = x^{(k)} - x^*$, rather than directly with the sequence of approximate solutions.  It is easy to see that the residual solutions satisfy a homogeneous version of the recurrence which defines the approximate solutions, i.e. $Mr^{(k+1)} = (M-A)r^{(k)}.$  Indeed, if $M$ is invertible, then the sequence of residuals is given by $r^{(k)} = C^kr^{(0)}$ where $C = I - M^{-1}A$ is known as the iteration matrix.  As a consequence, the iterative approaches to solving $Ax = b$ can be understood in terms of repeated matrix multiplication.  Indeed, an iterative method converges independent of the initial approximate solution if and of if the spectral radius of $C$, $\rho(C)$, is strictly less than 1.  Even more is true, specifically, if $\rho(C) < 1$ then $-\log_{10}(\rho(C))$ gives an lower bound on the rate of increase (in terms of number of iterations) of the approximate solutions, i.e. the number decimals of precision of $x^{(k)}$ is approximately $-\log_{10}(\rho(C)) k + \log_{10}(\norm{r^{(0)}})$. 

While iterative methods are quite effective and practical for solving linear systems, there are significant challenges that arise when considering large scale systems $Ax = b$ such as those that arise in many applications, such as the modeling of wind turbines, the modeling of reactor processes, computational chemistry for catalysis design, and many others. In particular, the matrix $A$ can be so large that it can not be effectively stored in memory (RAM) on a single compute node -- this necessitates either  taking a significant performance hit by repeatedly transferring portions of $A$ in and out of memory or dividing the system into parts and using many compute nodes to build the solution.  

More concretely, suppose that $A \in \real^{n \times n}$ and $b \in \real^{n}$ and we wish to use compute nodes labeled $1,2,\ldots,\ell$ to solve the system $Ax = b$.  One approach is to partition the rows of $A$ into $\ell$ sets resulting in rectangular matrices $A_i \in \real^{k_i \times n}$ and have a matching partition for $M$, $x$ and $b$.  In cases such as Jacobi iteration it is easy to extend the iterative scheme by considering the same partition on the diagonal matrix $M$, i.e. $M_ix_i^{(k+1)} = (M_i - A_i)x^{(k)} + b_i$, where $x^{(k)}$ is the approximate solution at iteration $k$ reconstructed from the local solutions $x_1^{(k)}, \ldots, x_{\ell}^{(k)}$. For non-diagonal $M$, such as with Gauss-Seidel iteration, the na\"ive approach would be to iteratively solve for the $x^{(k+1)}_i$ using the lower-triangular form and then ``pass" the information on $x^{(k+1)}_i$ onto the subsequent compute nodes so that the lower-triangular form can be used to solve $x_{i+1}^{(k+1)}, x_{i+2}^{(k+1)},\ldots,x_{\ell}^{(k+1)}.$  However, both of these approaches result in computational delays resulting from the time need to communicate between different compute nodes.  

The problem is especially acute in large high-performance computing (HPC) systems such as Aurora, Fugaku, and LUMI, where size of the system and over all network congestion can cause significant inter-node communication delays.
Indeed, for some scientific computations it is estimated that over 50\% of the computational time is taken up by communication requirements~\cite{jain2016optimization}.  In other contexts it is known that by performing work \emph{asynchronously}, that is, without waiting for the results of computation or work on remote nodes, can significantly increase the performance of HPC systems \cite{Suetterlein:introspector,Suetterlein:AMT,Suetterlein:roofline}.  However, little is known about how iterative methods for solving linear systems are effected by asynchrony.  More generally, we wish to understand the \emph{price of asyncrhony}, that is, the effect of asynchronous updates on the convergence rate of iterative methods.


\mySub{Problem Statement}

As a first step towards understanding the price of asynchrony, we propose a simple model of asynchronous updates.  To begin, we assume that the iterative method can be viewed through the recurrence $x^{(k+1)} = M x^{(k)}$ where $M$ is a $n \times n$ real matrix with $\rho(M) < 1$.  To build a simple model of asynchronous updates for $M$, we partition the matrix $M$ into an $\ell \times \ell$ array of blocks and for each $i \in [\ell]$ define the update scheme 
\[ \tilde{x}_i^{(k+1)} = \sum_j M_{ij}\tilde{x}_j^{(k - \delta_{ji})}, \] where $\delta_{ji}$ represents the ``delay" in information from compute node $j$ reaching compute node $i$ (and so we assume that $\delta_{ii}=0$).  The \emph{price of asynchrony} can be viewed as the difference between the rates at which $\norm{x^{(k)}}$ and $\norm{\tilde{x}^{(k)}}$ converge to 0.  With this framing, for any $\kappa \geq \max_{i,j} \delta_{ij}$ the asychronous updates can be viewed as a linear operator on $\real^{(\kappa+1)n}$ defined by  $(\tilde{x}^{(\kappa)}, \ldots, \tilde{x}^{(0)})$ is mapped to $(\tilde{x}^{(\kappa+1)}, \ldots, \tilde{x}^{(1)})$ where $\tilde{x}_i^{(\kappa+1)} = \sum_j M_{ij}\tilde{x}^{(\kappa - \delta_{ji})}_j$ for all blocks $i$.   We denote the matrix for this linear operator by $M^{(\delta,\kappa)}$ and note that $M^{(0,\kappa)}$ corresponds to the synchronous update which maintains a history of $\kappa$ steps.  

\begin{lemma}\label{L:spectra}
Let $M$ be a $n \times n$ real matrix with $\rho(M) < 1$ partitioned into an $\ell \times \ell$ array of blocks.  Let $\delta$ be an element of $\N^{\ell \times \ell}$ such that $\delta_{ii} = 0$ for all $i$.  Suppose that $\max_{ij} \delta_{ij} \leq \kappa \leq \kappa'$, then $\rho(M^{(\delta,\kappa)}) = \rho(M^{(\delta,\kappa')}).$ 
\end{lemma}

\begin{proof}
Suppose that $(\mu,v)$ is an eigenpair of $M^{(\delta,\kappa)}$.  Since $M^{(\delta,\kappa)}$ takes $(\tilde{x}^{(\kappa)}, \ldots, \tilde{x}^{(0)})$ to $(\tilde{x}^{(\kappa+1)}, \ldots, \tilde{x}^{(1)})$,
there is some vector $v_0 \in \real^n$ such that $v = (\mu^{\kappa}v_0, \mu^{\kappa-1}v_0, \ldots, \mu v_0, v_0).$  Further, since $v$ is an eigenvector of $M^{(\delta,\kappa)}$ we have that $\mu^{\kappa+1}(v_0)_i = \sum_j \mu^{\kappa - \delta^{ji}}v_0$.
Now define the vector $v' = (\mu^{\kappa}v_0, \mu^{\kappa - 1}v_0, \ldots, \mu v_0, v_0, \mu^{-1}v_0,\ldots, \mu^{\kappa - \kappa'} v_0)$ and consider $M^{(\delta,\kappa')}v'.$  By construction, we have that construction $M^{(\delta,\kappa')}v' =  (w', \mu^{\kappa}v_0,\ldots, \mu^{\kappa -\kappa'+1}v_0)$ for some $w' \in \real^n$.  Furthermore, $w'_i = \sum_j \mu^{\kappa - \delta_{ji}}v_0 = \mu^{\kappa+1}(v_0)_i$.  As a consequence $(\mu,v')$ is a eigenpair for $M^{(\delta,\kappa')}$. A similar argument shows that any eigenpair for $M^{(\delta,\kappa')}$ has a corresponding eigenpair for $M^{(\delta,\kappa)}$. In particular, the two matrices are co-spectral and hence have the same spectral radius.
\end{proof}

As an immediate consequence of Lemma \ref{L:spectra}, we have that $\rho(M^{(0,\kappa)}) = \rho(M)$ for any $\kappa.$
Furthermore, it suffices to restrict our attention to the spectral radius of $M^{(\delta,\kappa)}$ for any system of delays $\delta$ and $\kappa = \max \delta$. Thus, we will abuse notation slightly denote by $M^{(\delta)}$ any member of the family $\set{ M^{(\delta, k)} \mid k \geq \max \delta}$. 

Intuitively, it seems that if we fix a maximum delay $k$, then $\rho(M^{(\delta)})$ is maximized when the delay between any pair of nodes achieves the maximum delay.  In this case, uniformity of the delays allows for an explicit calculation of $\rho(M^{(\delta)})$.

\begin{lemma}\label{L:uniform}
Let $M$ be an $n \times n$ real matrix with $\rho(M) < 1$ partitioned into an $\ell \times \ell$ array of blocks. Let $\delta$ be an element of $\N^{\ell \times \ell}$ such that $\delta_{ij} = k$ for all non-zero blocks of $M$, then $\rho(M^{(\delta)}) = \rho(M)^{\nicefrac{1}{k+1}}$.
\end{lemma}

\begin{proof}
    As noted in Lemma \ref{L:spectra}, any eigenpair $(\mu,v)$ of $M^{(\delta,\kappa)}$ has the form $v = (w, \mu^{-1},\ldots,\mu^{-\kappa}w)$, where $\mu w_i = \sum_j \mu^{-\delta_{ij}}M_{ij} w_j$.  However, since $\delta_{ij} = k$ for any non-zero block of $M$, this implies that $\mu w = \mu^{-k} M w$ and $(\mu^{k+1},w)$ is an eigenpair for $M$.  It is easy to see that a similar construction will take an eigenpair for $M$ to an eigenpair for $M^{(\delta,\kappa)}$.  The desired result on the spectral radius follows immediately.
\end{proof}

This leads naturally to the following open questions:
\begin{problem}
Let $M$ be a $n \times n$ real matrix with $p(M)= \lambda < 1$ partitioned into an $\ell \times \ell$ array of blocks.  Is it the case that for all elements $\delta$ of $[k]^{\ell \times \ell}$ with zero diagonals, $\rho(M^{(\delta)}) \leq \lambda^{\nicefrac{1}{k+1}}?$  If not, is there some function $f \colon (0,1) \times \N \rightarrow (0,1)$ such that $\rho(M^{(\delta)}) \leq f(\lambda,k)$?
\end{problem}
Intuitively, it seems as if this is the worst possible case for the convergence rate, that is, if the all the delays are at most $k$. 
\mySub{Experimental Results}

\begin{figure}
    \centering
    \hfill
        \subfloat[Single Delays\label{F:single}]{
    \begin{tikzpicture}[decoration={
    markings,
    mark=at position 0.5 with {\arrow{>}}}]
    \node[circle,draw] (4) at (-2,-1.155) {4};
    \node[circle,draw] (3) at (0,2.31) {3};
    \node[circle,draw] (2) at (2,-1.155) {2};
    \node[circle,draw] (1) at (0,0)  {1};
    \draw[postaction = {decorate}] (1) to [bend right=10] node[midway,left] {\tiny 5} (2);
    \draw[postaction = {decorate}] (1) to [bend right=10] node[midway,below right] {\tiny 5} (3);
    \draw[postaction = {decorate}] (1) to [bend right=10] node[midway,above] {\tiny 5} (4);
    \draw[white,postaction = {decorate}] (4) to [bend right=10] node[midway,below] {\tiny 1} (2);
    \end{tikzpicture}
    }
    \hfill
        \subfloat[Poisson Delays\label{F:poisson}]{
    \mbox{
    \begin{tikzpicture}[decoration={
    markings,
    mark=at position 0.5 with {\arrow{>}}}]
    \node[circle,draw] (4) at (-2,-1.155) {4};
    \node[circle,draw] (3) at (0,2.31) {3};
    \node[circle,draw] (2) at (2,-1.155) {2};
    \node[circle,draw] (1) at (0,0)  {1};
    \draw[postaction = {decorate}] (1) to [bend right=10] node[midway,left] {\tiny 1} (2);
    \draw[postaction = {decorate}] (1) to [bend right=10] node[midway,below right] {\tiny 3} (3);
    \draw[postaction = {decorate}] (1) to [bend right=10] node[midway,above] {\tiny 2} (4);
    \draw[postaction = {decorate}] (2) to [bend right=10] node[midway,above] {\tiny 4} (1);
    \draw[postaction = {decorate}] (2) to [bend right=10] node[midway,above right] {\tiny 2} (3);
    \draw[postaction = {decorate}] (2) to [bend right=10] node[midway,above] {\tiny 5} (4);
    \draw[postaction = {decorate}] (3) to [bend right=10] node[midway,below left] {\tiny 2} (1);
    \draw[postaction = {decorate}] (3) to [bend right=10] node[midway,below] {\tiny 1} (2);
    \draw[postaction = {decorate}] (3) to [bend right=10] node[midway,above left] {\tiny 2} (4);
    \draw[postaction = {decorate}] (4) to [bend right=10] node[midway,right] {\tiny 4} (1);
    \draw[postaction = {decorate}] (4) to [bend right=10] node[midway,below] {\tiny 2} (2);
    \draw[postaction = {decorate}] (4) to [bend right=10] node[midway,below] {\tiny 5} (3);
    \end{tikzpicture}
    }
    }
    \hfill\phantom{}
    \caption{Delay patterns for asynchronous computation.  A weight $d$ edge directed from $j$ to $i$ represents a delay of $d$ time steps in information propagating from $j$ to $i$, i.e. $\delta_{ji} = d$.}\label{F}
\end{figure}
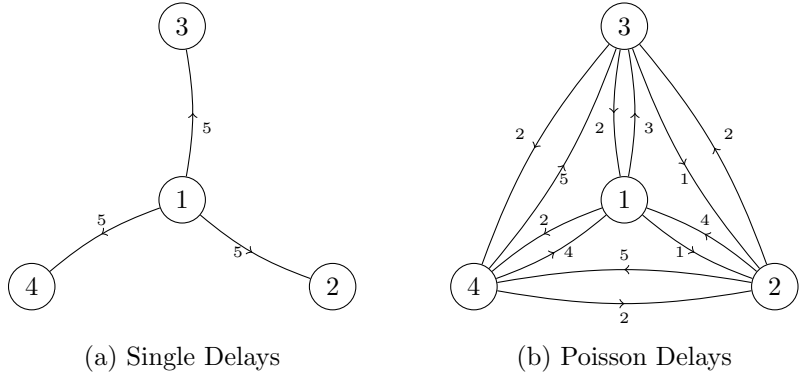

\begin{figure}[h!]
    \centering
    \hfill
    \subfloat[Single Delay\label{F:single_plot}]{\includegraphics[trim = 25 15 45 35, clip, width=.40\textwidth]{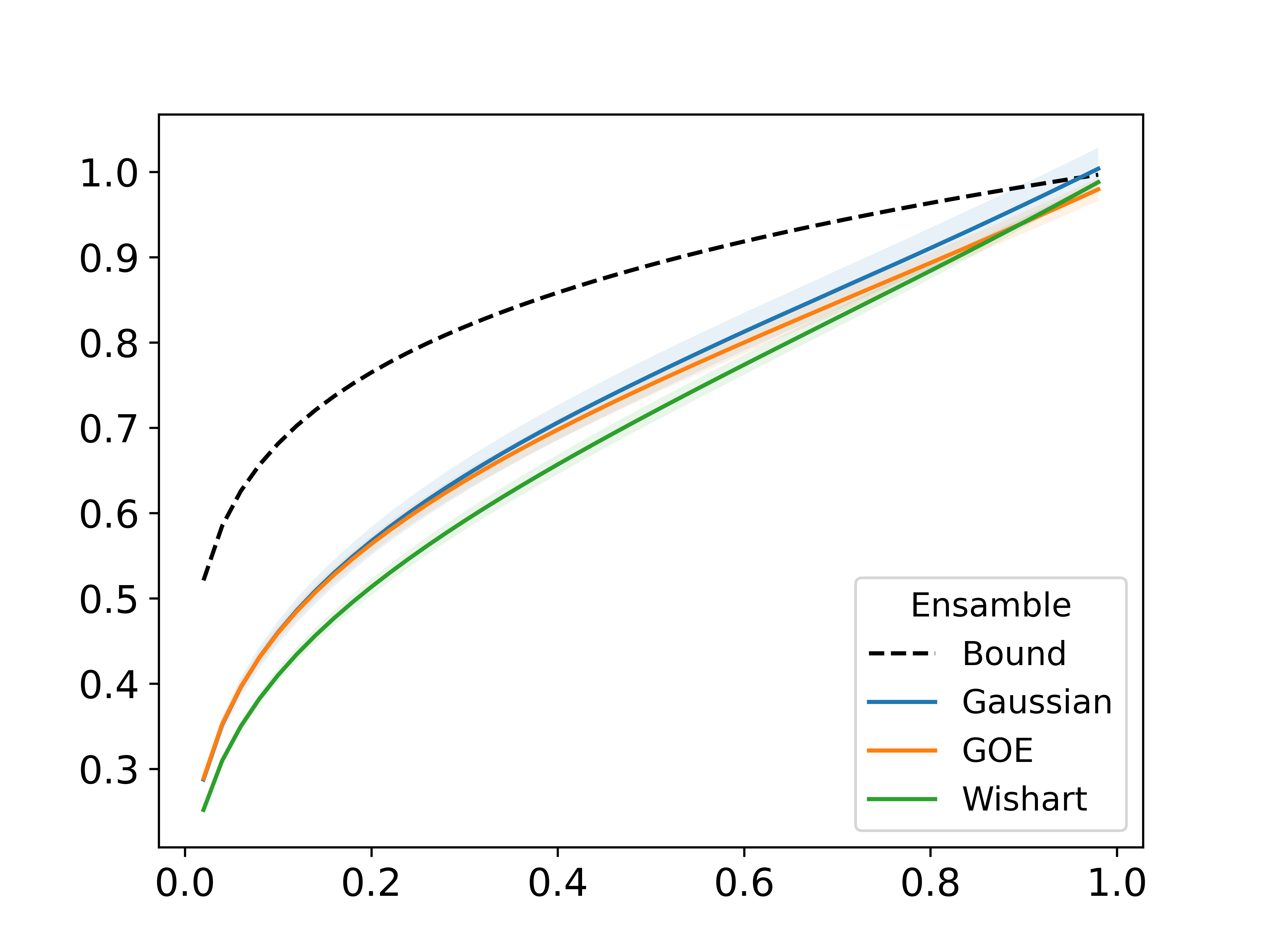}} \hfill
    \subfloat[Poisson Delays\label{F:poisson_plot}]{\includegraphics[trim = 25 15 45 35, clip, width=.40\textwidth]{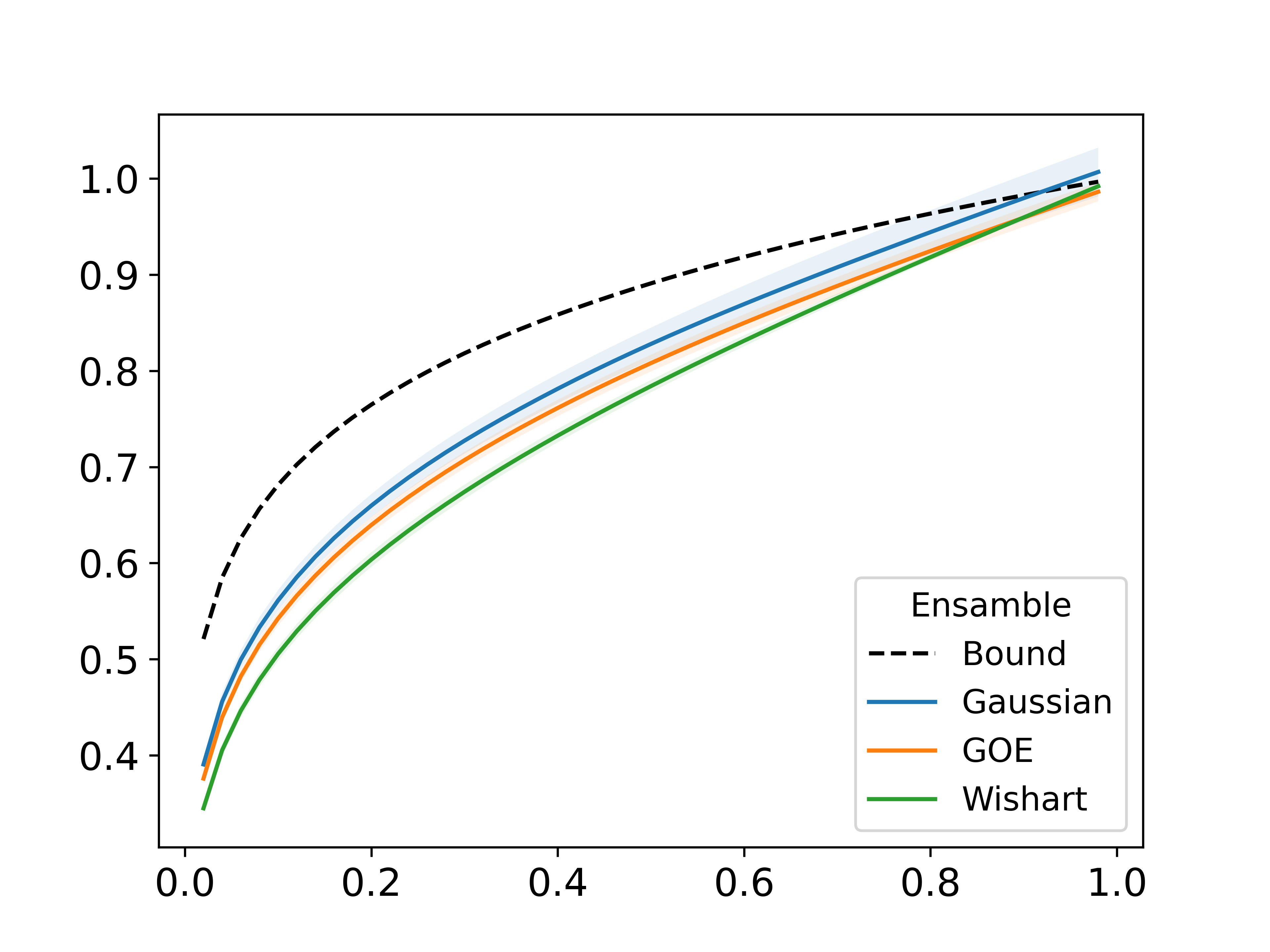}} \hfill \phantom{} \\
    \hfill
    \subfloat[Single Delay, Block Jacobi\label{F:jac_single_plot}]{\includegraphics[trim = 25 15 45 35, clip, width=.40\textwidth]{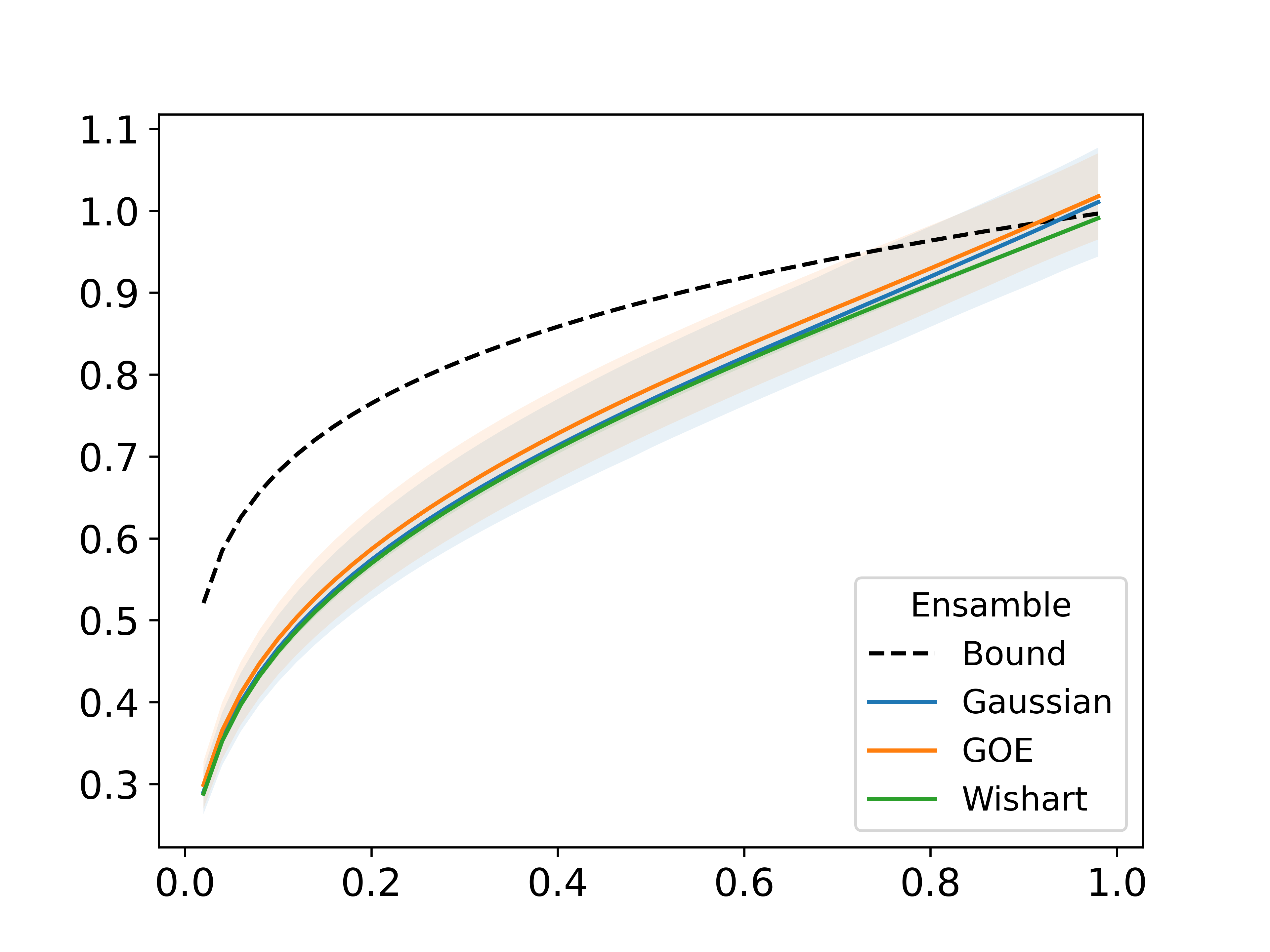}} \hfill
    \subfloat[Poisson Delays, Block Jacobi\label{F:jac_poissoin_plot}]{\includegraphics[trim = 25 15 45 35, clip, width=.40\textwidth]{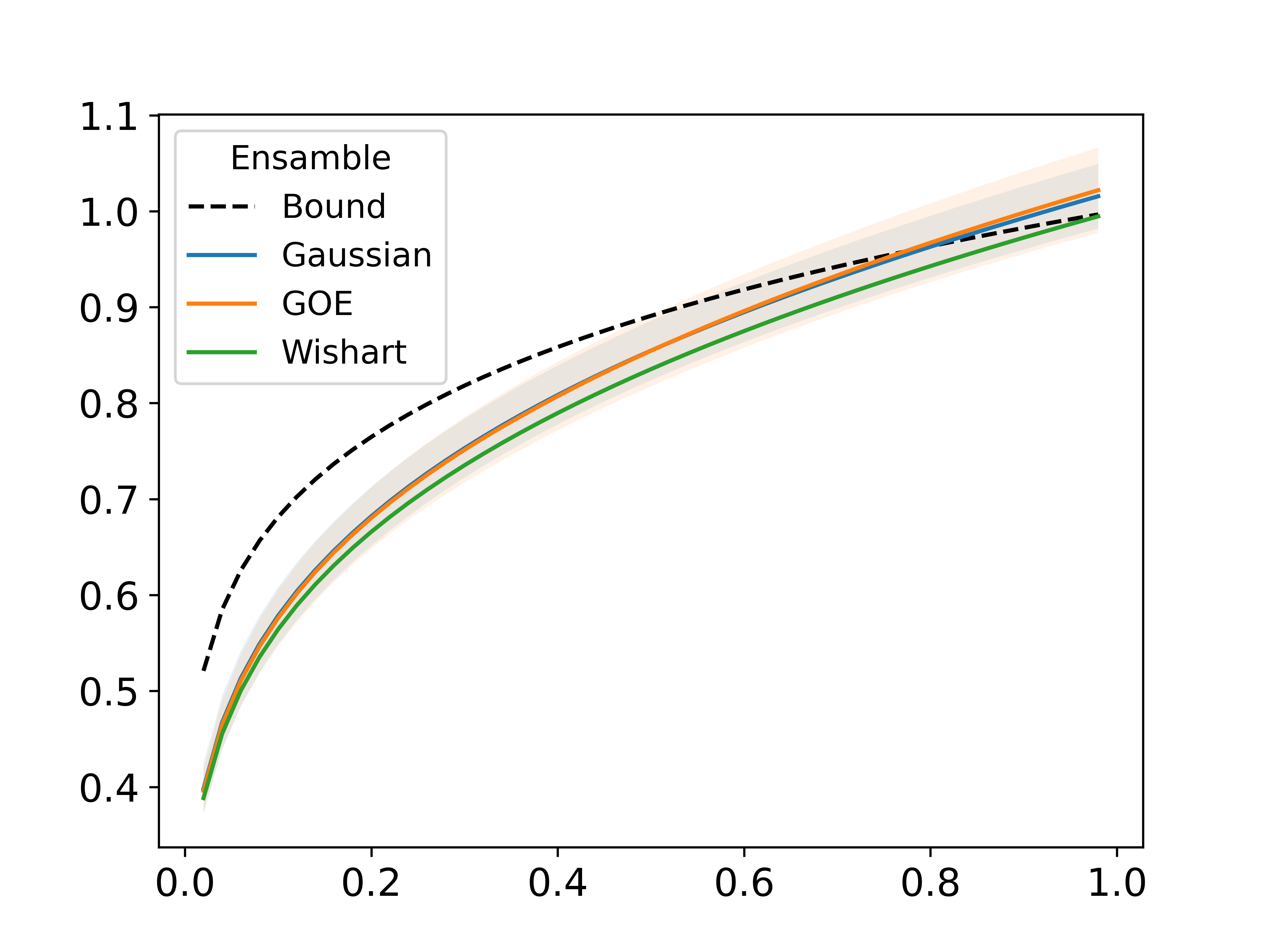}} \hfill \phantom{}
    \caption{{Price of Asynchrony.}
    For each of the three matrix ensembles, the solid line represents the average spectral radius over 50 different samples from the ensemble, while the shaded area represents a range of two standard deviations.  The dashed line is the conjectured upper bound in terms of the spectral radius of the synchronous update matrix.}
    \label{F:spectral_plots}
\end{figure}

As a first pass towards resolving this problem, consider a random matrix $M$ selected from one of three different ensembles derived from a matrix $X$ with independent normally distributed entries with mean zero;  $M = X$,  the Gaussian Orthonormal Ensemble where $M = X + X^T$, and the Wishart ensemble where $M = XX^T$.  For notational convenience, we assume that $M$ is re-scaled to having unit spectral radius.   We also consider, for each of these matrix ensembles, the iteration matrix associated with a block Jacobi iteration (again, normalized to have norm 1), that is, the diagonal blocks are zero and the off-diagonal blocks are given by $M_{ii}^{-1}M_{ij}$ where $M_{xy}$ denotes the $x,y$ block in $M$.  In each of these six, cases we also consider two different different delay patterns with 4 compute nodes;  a single node with a delay of 5 and all others having a delay of 0 (Figure \ref{F:single}) and randomly generated Poisson delays with a mean of 3 (Figure \ref{F:poisson}).  In Figure \ref{F:spectral_plots}, we plot the relationship between $\rho((cM)^{(\delta)})$ versus $c \in (0,1)$ for each of these graphs over 50 different trials. 

These experiments would seem to provide at least some support for the conjecture that $\rho(M^{(\delta)}) \leq \rho(M)^{\nicefrac{1}{\max \delta}},$ at least when $\rho(M)$ is sufficiently small.  In the case where $\rho(M)$ is close to one, we expect that at least part of the violation of the conjecture may be as a result of numerical instabilities in the spectral radius calculation.  In particular, from the proofs of Lemma \ref{L:spectra} and \ref{L:uniform}, it easy to see that despite $M^{(\delta,\kappa)}$ being an $(\kappa +1)n \times (\kappa+1)n$ matrix, the sum of the geometric multiplicities is at most $n$.  To further complicate matters, $M^{(\delta)}$ is non-Hermitian and has entries of both positive and negative signs, increasing the likelihood that the eigenvalue which yields the spectral radius is complex.

It is also interesting to note that there is relatively little difference between the behavior of the full matrix iteration and the block Jacobi iteration when there is a single constant delay, while for Poisson delays the block Jacobi iteration has a significantly higher spectral radius.  One possible explanation is to consider the average effective delays (define as the mean delay over non-zero blocks of $M$) which is $\nicefrac{15}{16}$ for a single delayed node, $\nicefrac{15}{12}$ for Jacobi iteration with a single delayed node, $\nicefrac{33}{16}$ for Poisson delays, and $\nicefrac{33}{12}$ for Jacobi iteration with Poisson delays.  As a result, one might speculate that the important parameter for understanding $\rho(M^{(\delta)})$ is the mean delay, perhaps weighted by the spectral radius of the individual blocks.  

\mySub{Acknowledgements} The authors gratefully acknowledge the funding support from the Applied Mathematics Program within the U.S. Department of Energy’s Office of Advanced Scientific Computing Research as part of RAndomized Techniques For Iterative Solvers in Heterogeneous environments (RATFISH). Pacific Northwest National Laboratory is operated by Battelle for the DOE under Contract DE-AC05-76RL01830.


\printbibliography[heading=subbibliography]
\end{refsection}

\begin{refsection}

\newcommand{\LZF}[2][\ell]{Z_{(#1)}\!\paren{#2}}

\section{$(n-k)$-contingent Zero-Forcing for Power Grids}\label{sec:aksoy}
\begin{flushright}
{\it Sinan G.\ Aksoy, Anthony V.\ Petyuk, Sandip Roy, Stephen J.\ Young}
\end{flushright}

\mySub{Introduction}
With the increasing penetration of microgrids, renewable energy sources, power-generation at the edge, and distributed energy resources, maintaining the stability of power grid transmission has changed dramatically.  For example, the majority of grid power generation used to be achieved by spinning large turbines (e.g. hydro-electric turbines spun by falling water, steam powered turbines in gas/coal/nuclear power plants).  However, in 2022 almost 14\% of the utility scale generation was from wind or solar~\cite{electric_power}. In contrast to turbine-based generation, neither wind nor solar power has any ``inertia" and so is subject to rapid changes in the total generation which are not controlled by operators.  This puts significant pressure on control schemes designed to manage the stability of the power grid through the regulation of power generation.  Such control schemes are also made challenging by the increasing penetration of microgrids in the power system: these are small, regional areas of the power grid which contain sufficient internal generation that disconnect themselves from the broader power grid in order to improve local stability.  As microgrids change their connectivity to the power grid, the fundamental electrical equations which govern power flow on the grid also change, impacting control schemes.

The changing grid not only poses challenges for stable control, but also presents several opportunities to increase both resilience and efficiency.  For example, phasor measurement units (PMUs) have been deployed at selected places within the power-grid over the last 25 years, providing measures of the local power grid state (i.e., voltage magnitude and phase angle) up to 120 times a second -- a significant increase over previous methods yielding one measurement every few seconds.  
Another opportunity to increase the resilience and stability of the grid is afforded by the increasing penetration of distributed energy resources (DERs), such as grid-scale battery installations.
Not only do these large battery installations provide a means to temporally arbitrage power from renewable sources (e.g., store excess power generated during the day by solar panels to use at night), they also provide alternative points of control for the power grid via selective charging, islanding, and discharging to rapidly inject power from the overall grid system.  
Furthermore, because of the significantly smaller footprint of grid-scale batteries (as compared to traditional power plants) it is possible to deploy these resources in locations which are advantageous to control and stability. 

Given this variety of power-grid structures and operating points, it might seem nearly intractable to decide whether a collections of PMUs provide sufficient observability of the state of grid, or whether a collection DERs are able to provide stabilization for the grid.  Fortunately, both of these problems can be reduced to a combinatorial problem on the structure of the grid known as zero-forcing.

\begin{definition}[Zero-Forcing Move]
Given a graph $G = (V,E)$ and a set $S \subseteq V$ of vertices colored blue.  A vertex $v \in V-S$ can be colored by a \emph{zero-forcing move} if there is a vertex $s$ such that $\set{v}  = N(s) - S$, that is, if there is a vertex $s \in S$ such that  $v$ is the unique neighbor of $s$ not colored blue. If $v$ is colored via a zero-forcing move that relies on $s$, we will say that $s$ is used to color $v$.
\end{definition}

\begin{definition}[Zero-Forcing Set]
Given a graph $G = (V,E)$, a set $Z \subseteq V$ is called a \emph{zero-forcing set} if there is a sequence of zero-forcing moves such that all vertices in $V$ are colored.
\end{definition}

There is an extensive literature around determining the size and structure of minimal zero-forcing sets (see for instance \cite{ZFbook} and the references therein). Concerning power-grid applications, \cite{monshizadeh2014zero} shows the linear dynamics associated to the power-flow equations are controllable by controllers at a zero-forcing set, and \cite{Brueni:PMUplacement,smith2020optimal} shows if phasor-measurement units are placed at a power-dominating set\footnote{A power-dominating set is a set of vertices $S$ such that, combined with their neighborhood, they form a zero-forcing set.} then the phase of every bus in the network is recoverable. 
Consequently, distributed energy resource placement (such as grid-scale batteries) and phasor-measurement units may enable the real-time steering of the power-grid to increase efficiency and resiliency. In addition, recent work~\cite{Roy:zero_forcing} has shown zero-forcing sets can be used to identify control points for microgrids.  However, these efforts fail to take into account the dynamic nature of the grid.  In particular, zero-forcing is defined for a static graph which does not allow for changes in the underlying network topology, limiting its applicability to the observability and control of the power grid.

\mySub{Contingent Zero-Forcing}

To address the limitations of zero-forcing in the context of the power grid, 
we propose the following extension of the zero-forcing to identify resilient sets which can achieve the observability and control of the grid.
\begin{definition}[$(n-k)$-contingent Zero Forcing]
Given a graph $G = (V,E)$, a set $Z \subseteq V$ is a $(n-k)$-contingent zero-forcing set if for every set of edges $E'$ of size at most $k$, $Z$ is a zero-forcing set for $G' = (V, E-E')$.
\end{definition}

We note the edge-leaky zero-forcing number introduced independently by Alameda, Kritschgau, and Young~\cite{MYoung:EdgeLeakyZF} is essentially equivalent to this definition\footnote{In contrast to $(n-k)$-contingent zero-forcing, edges in the edge-leaky zero-forcing process aren't removed from the graph, rather their use is prohibited. Consequently, it is possible that in some stage of the zero-forcing process the contingent zero-forcing process can force a strict superset of the vertices forcible in the equivalent edge-leaky process.  However, Alameda, Kritschgau, and Young go on to show the $\ell$-edge-leaky zero-forcing sets are the same as $\ell$-leaky zero-forcing sets. Thus, by Lemma \ref{L:equiv} $\ell$-edge leaky zero-forcing is equivalent to $(n-\ell)$-contingent zero-forcing.}. Since the transmission and distribution level of power grids frequently exhibit ``tree-like" structure, it is natural to first consider nature of the $(n-k)$-contingent zero-forcing sets on trees.  To explore this, we first recall the following fact about the zero-forcing sets of trees.

\begin{theorem}[\cite{Kenter:LeakyZF}]\label{T:treeZF}
If $T$ is a tree with $t \geq 2$ leaves, then any collection of at least $t-1$ leaves is a zero-forcing set for $T.$
\end{theorem}

In the setting of $(n-k)$-contingent zero-forcing, we apply this to prove the following lemma. 

\begin{lemma}\label{L:tree_kZF}
    Let $T = (V,E)$ be a tree. For any $k \geq 1$, a set $Z \subseteq V$ is a $(n-k)$-contingent zero-forcing set if and only if $Z$ contains all vertices of degree at most $k$.  
\end{lemma}

\begin{proof}
    If $Z$ is a $(n-k)$-contingent zero-forcing set for any graph, then it must contain all vertices of degree at most $k$.  In particular, if any vertex $v$ of degree at most $k$ is not in $Z$ then removing all of the incident edges to $v$ results in a graph where $Z$ is not a zero-forcing set.  

    Now let $Z$ be the set of vertices of degree at most $k$ in $T$ and let $F$ be the forest formed by removing a set of $k$ edges from $T$.  We first note that any isolated vertex in $F$ has degree at most $k$ in $T$ as we removed at most $k$ edges. Suppose there is some tree $T'$ which has two leaves $u,v \not\in Z$. We note that since $u$ and $v$ are in the same tree in $F$, there edge $\set{u,v}$ is not in $T$.  Thus, in order for both $u$ and $v$ to be leaves in $F$ at least $k + k > k$ edges must be removed, contradicting the construction of $F$.  Thus, every tree in $F$ is either an isolated vertex belong to $Z$ or a tree with at most one leaf not in $Z$. Thus,  by Theorem \ref{T:treeZF}, every tree in $F$ is forcible by $Z$.  
\end{proof}

The idea of $(n-k)$-contingent zero-forcing is the not the first attempt to develop a notion of zero-forcing that is resilient to changes in the underling graph. For example, Dillman and Kenter~\cite{Kenter:LeakyZF} consider the resilience in the context of water flows where water can ``leak" out at the joints preventing vertices from participating in the forcing process.  More formally, they give the following definition:

\begin{definition}[$\ell$-leaky Zero-Forcing]
Given a graph $G = (V,E)$, a set $Z \subseteq V$ is an $\ell$-leaky zero-forcing set if for every set of vertices $L \subseteq V$ of size at most $\ell$, $Z$ is a zero-forcing set of $G'$, where $G'$ is formed by adding a pendant vertex to each vertex in $L$.
\end{definition}

Observe that any vertex $v \in L$ can not be used to color any other vertex via a zero-forcing move. In other words, the zero-forcing process ``leaks" out of the graph on the vertices in $L$ so they are only able to result in the coloring of new pendant vertices.  In many ways, the leaky zero-forcing process can be thought of as an vertex version of the contingent zero-forcing process: in leaky zero-forcing an unknown set of vertices can not participate in the zero-forcing process, while in the contingent zero-forcing process an unknown set of edges can not participate in the zero-forcing process. With this observation, it is unsurprising these two processes are effectively equivalent, as we show below.

\begin{lemma}\label{L:equiv}
A set $S$ is a $(n-k)$-contingent zero-forcing set if and only if it is an  $k$-leaky zero-forcing set.
\end{lemma}

\begin{proof}
Let $G$ be a graph and suppose $Z$ is a $(n-k)$-contingent zero-forcing set for $G.$ Fix an arbitrary set of $L$ leaks in the graph $G$, where $\size{L} \leq k$.  We will identify a set $\mathcal{E}$ of at most $\ell$ edges such that there is a forcing sequence of $G - \mathcal{E}$ which does not use any vertex in $L$.  To this end, we iteratively build the set $\mathcal{E}$ of removed edges based on $Z$.  Note that as long as $\size{\mathcal{E}} \leq k$, the forcing process from $Z$ can continue as $Z$ is a $(n-k)$-contingent zero-forcing set.   Now, to identify these edges of $\mathcal{E}$, consider sequentially applying the forcing process until the first time a vertex $v \in L$ is used to color $u$.  At this point, add the edge $\set{v,u}$ to $\mathcal{E}$ and remove it from the graph $G$.  Since all of the neighbors of $v$ are now colored, there is no future step in the forcing process in which $v$ colors a neighbor.  Thus we may effectively remove $v$ from $L$.  Repeating this process results in a set $\mathcal{E}$ of size at most $k$ and a sequence of zero-forcing moves for $G - \mathcal{E}$, such that no vertex in $L$ forces any other vertex. Thus $S$ is a forcing set for $G$ with leaks at all the vertices in $L$.  

The converse proceeds similarly.  Let $Z$ be a $k$-leaky zero-forcing set for $G$ and fix an arbitrary set of $\mathcal{E}$ edges in $G$.   We again interatively apply the zero-forcing process and build a set of leaks $L$ as needed.  Specifically, for each edge $\set{u,v} \in E$ we add the first vertex from the edge that is colored to the set of leaks $L$.  Clearly, this ensures that the edge $\set{u,v}$ is not used in the forcing process and adds at most $k$ vertices to the set $L$. 
\end{proof}

In general, determining the size of the minimal zero-forcing set is $\mathcal{NP}$-\texttt{complete}, even for highly-restricted classes of graphs (see \cite{Fallat:PSD_ZeroForcing}, for instance).  Thus, it is exceedingly likely that determining the size of the minimal $(n-k)$-contingent zero-forcing number is also $\mathcal{NP}$-\textrm{complete}.  However, the class of graphs corresponding to power-grid networks is known to have a number of unusual structural properties (see, for instance \cite{Aksoy2018,Young:NoN}), including average degree between 1 and 2, diameter scaling like the square root of the number of vertices, and the presence of long-cycles.  In other words, typical power-grid networks are tree-like by Lemma \ref{L:tree_kZF} the minimal $(n-k)$-contingent zero-forcing sets are known exactly for trees.  This leads naturally to our first two open problems:

\begin{problem}
    Is there a structural characterization of the minimal-sized $(n-k)$-contingent zero-forcing sets for graphs with few cycles? For sparse graphs with large diameter?
\end{problem}

\begin{problem}
    Is there an efficient algorithm to determine the minimal-sized $(n-k)$-contingent zero-forcing set for graphs with few cycles? For sparse graphs with large diameter?
\end{problem}

Even if the structural properties of power-grid networks are insufficient to yield an efficient means of computing a minimal, or nearly so, $(n-k)$-contingent zero-forcing set, if such a set is sufficiently small the additional computational effort may be worth while.  However, in order to assess the value of identifying the minimal $(n-k)$-contingent zero-forcing sets, it would be helpful to have a rough estimate of size of the resulting set.  

\begin{problem}
    What is the typical size, in terms of the number of vertices and edges, of the minimal $(n-k)$-contingent zero-forcing sets in graphs exhibiting structural properties characteristic of power grids?
\end{problem}

While this question may be answerable by restricting attention to the class of graphs with a few hallmark structural properties (e.g., limitation on density, diameter, maximum degree, etc.) it may be more useful to build off the recent work characterizing the zero-forcing number for the Erd\H{o}s-Reny\'i random graph~\cite{Pralat:ZeroForcingRandom} and random regular graphs~\cite{Pralat:ZF_regular,Sudakov:ZeroForcing}.  While the structural properties of these random graph models are far from those of the power-grid, there are a few recent models such as the Chung-Lu Chain~\cite{Aksoy2018} and the Geometric Delaunay~\cite{Young:GeoDe} models which have been shown to capture an array of structural properties of the power grid.

\mySub{Acknowledgements}
Sinan G.\ Aksoy, Anthony V.\ Petyuk, and Stephen J.\ Young were supported by the Resilience through Data-driven Intelligently-Designed Control (RD2C) Initiative, under the Laboratory Directed Research and Development (LDRD) Program at Pacific Northwest National Laboratory (PNNL).  PNNL is a multi-program national laboratory operated for the U.S. Department of Energy (DOE) by Battelle Memorial Institute under Contract No. DE-AC05-76RL01830. Sandip Roy conducted this research during an appointment at the U.S. National Science Foundation, supported by an Intergovermental Personnel Act agreement with Washington State University.
\printbibliography[heading=subbibliography]
\end{refsection}



\end{document}